\documentclass[format=acmsmall, screen=true, review=false, authorversion=true]{acmart}
\usepackage{mystyle, hyperref, multirow}

\makeatletter                   
\def\mdseries@tt{m}             
\makeatother                    

\usepackage[draft=true]{minted}
\usepackage{xcolor,cancel}

\definecolor{lightgray}{gray}{.95}
\newminted{python}{
    linenos=true,
    bgcolor=lightgray,
    tabsize=4,
    fontfamily=courier,
    fontsize=\small,
    xleftmargin=5pt,
    xrightmargin=5pt
}

\usepackage{hyperref}           
\hypersetup{
    colorlinks=true,
    linkcolor=blue,
    filecolor=red,      
    urlcolor=magenta,
    breaklinks=true,            
}
\usepackage{breakurl}           

\let\vec\relax
\DeclareMathOperator*{\vec}{\textbf{vec}}
\DeclareMathOperator*{\epi}{epi}
\DeclareMathOperator*{\dx}{dx}
\DeclareMathOperator*{\dt}{dt}
\DeclareMathOperator*{\dS}{dS}
\DeclareMathOperator*{\curl}{curl}
\DeclareMathOperator*{\argmin}{argmin}
\graphicspath{{./pic/},{./results/}}
\begin{document}
\sloppy

\title[Automating the formulation and resolution of convex variational problems]{Automating the formulation and resolution of convex variational problems:  applications from image processing to computational mechanics}
\author{Jeremy Bleyer}
\email{jeremy.bleyer@enpc.fr}
\orcid{0001-8212-9921}
\affiliation{%
  \institution{Laboratoire Navier UMR 8205 (ENPC-IFSTTAR-CNRS), Universit\'e Paris-Est}
  \streetaddress{6-8 av Blaise Pascal, Cit\'e Descartes}
  \city{Champs-sur-Marne}
  \country{France}
  \postcode{77455}
}


\begin{abstract}
Convex variational problems arise in many fields ranging from image processing to fluid and solid mechanics communities. Interesting applications usually involve non-smooth terms which require well-designed optimization algorithms for their resolution. The present manuscript presents the Python package called \texttt{fenics\_optim} built on top of the FEniCS finite element software which enables to automate the formulation and resolution of various convex variational problems. Formulating such a problem relies on FEniCS domain-specific language and the representation of convex functions, in particular non-smooth ones, in the conic programming framework. The discrete formulation of the corresponding optimization problems hinges on the finite element discretization capabilities offered by FEniCS while their numerical resolution is carried out by the interior-point solver Mosek. Through various illustrative examples, we show that convex optimization problems can be formulated using only a few lines of code, discretized in a very simple manner and solved extremely efficiently. 
\end{abstract}


\maketitle

\section{Introduction}

Convex variational problems represent an important class of mathematical abstractions which can be used to model various physical systems or provide a natural way of formulating interesting problems in different areas of applied mathematics. Moreover, they also often arise as a relaxation of more complicated non-convex problems. Optimality conditions of constrained convex variational problems correspond to variational inequalities which have been the topic of a large amount of work in terms of analysis or practical applications \cite{duvaut2012inequalities, kinderlehrer1980introduction, tremolieres2011numerical}.

In this work, we consider convex variational problems defined on a domain $\Omega \subset \RR^d$ ($d=2,3$ for typical applications) with convex constraints of the following kind:
\begin{equation}
\begin{array}{rl} \displaystyle{\inf_{u\in V}} & J(u) \\
\text{s.t.}  &u\in \Kk
\end{array} \label{co-prob}
\end{equation}
where $u$ belongs to a suitable functional space $V$, $J$ is a convex function and $\Kk$ a convex subset of $V$. Some variational inequality problems are formulated naturally in this framework such as the classical Signorini obstacle problem in which $\Kk$ encodes the linear inequality constraint that a membrane displacement cannot interpenetrate a fixed obstacle (see section \ref{sec:obstacle}). An important class of situations concerns the case where $J$ can be decomposed as the sum of a smooth and a non-smooth term. Such a situation arises in many variational models of image processing problems such as image denoising, inpainting, deconvolution, decomposition, etc. In some cases, such as limit analysis problems in mechanics for instance, smooth terms in $J$ are absent so that numerical resolution of \eqref{co-prob} becomes very challenging \cite{tremolieres2011numerical, kanno2011nonsmooth}. Important problems in applied mathematics such as optimal control \cite{malanowski1982convergence} or optimal transportation \cite{benamou2000computational,villani2003topics,papadakis2014optimal,peyre2019computational} can also be formulated, in some circumstances, as convex variational problems. This is also the case for some classes of topology optimization problems \cite{bendsoe_topology_2004}, which can also be extended to non-convex problems involving integer optimization variables \cite{kanno2010mixed,yonekura2010global}. Finally, robust optimization in which optimization is performed while taking into account uncertainty in the input data of \eqref{co-prob} has been developed in the last decade \cite{ben2002robust, ben2009robust}. It leads, in some cases, to tractable optimization problems fitting the same framework, possibly with more complex constraints. 

%


The main goal of this paper is to present a numerical toolbox for automating the formulation and the resolution of discrete approximations of \eqref{co-prob} using the finite-element method. The large variety of areas in which such problems arise makes us believe that there is a need for a versatile tool  which will aim at satisfying three important features:
\begin{itemize}
\item straightforward formulation of the problem, mimicking in particular the expression of the continuous functional;
\item automated finite-element discretization, supporting not only standard Lagrange finite elements but also DG formulations and $H(\div)/H(\curl)$ elements;
\item efficient solution procedure for all kinds of convex functionals, in particular non-smooth ones.
\end{itemize}

In our proposal, the first two points will rely extensively on the versatility and computational efficiency of the FEniCS open-source finite element library \cite{logg2012automated,AlnaesBlechta2015a}. FEniCS is now an established collection of components including the DOLFIN C++/Python Interface \cite{LoggWells2010a,LoggWellsEtAl2012a}, the Unified Form Language \cite{Alnaes2012a,AlnaesEtAl2012}, the FEniCS Form Compiler \cite{KirbyLogg2006a,LoggOlgaardEtAl2012a}, etc. Using the high-level DOLFIN interface, the user is able to write short pieces of code for automating the resolution of PDEs in an efficient manner. For all these reasons, we decided to develop a Python package called \texttt{fenics\_optim} as an add-on to the FEniCS library. We will also make use of Object Oriented Programming possibilities offered by Python for defining easily our problems (see \ref{sec:libary}). Our proposal can therefore be considered to be close to high-level optimisation libraries based on \textit{disciplined convex programming} such as CVX\footnote{\url{http://cvxr.com/cvx/}} for instance. However, here we really concentrate on the integration within an efficient finite-element library offering symbolic computation capabilities. As mentioned later, integration with other high-level optimisation libraries will be an interesting development perspective.\\

Concerning the last point on solution procedure, we will here rely on the state-of-the art conic programming solver named Mosek \cite{mosek}, which implements extremely efficient primal-dual interior-point algorithms \cite{andersen2003implementing}. Let us mention first that there is no ideal choice concerning solution algorithms which mainly depends on the desired level of accuracy, the size of the considered problem, the sparsity of the underlying linear operators, the type of convex functionals involved, etc. In particular, in the image processing community, first-order proximal algorithms are widely used since they work well in practice for large scale problems discretized on uniform grids \cite{papadakis2014optimal}. Moreover, high accuracy on the computed solution is usually not required since one aims mostly at achieving some decrease of the cost function but not necessarily an accurate computation of the optimal point. In contrast, interior-point methods can achieve a desired accuracy on the solution in polynomial time and, in practice, in quasi-linear time since the number of final iterations is often observed to be nearly independent on the problem size. However, as a second-order method, each iteration is costly since it requires to solve a Newton-like system. Iterative solvers are also difficult to use in this context due to the strong increase of the Newton system conditioning when approaching the solution. As result, such solvers usually rely on direct solvers for factorizing the resulting Newton system, therefore requiring important memory usage. Nevertheless, interior-point solvers are extremely robust and quite efficient even compared to first order methods in some cases. For these reasons, the present paper will not focus on comparing different solution procedures and we will use only the Mosek solver but including first-order algorithms in the \texttt{fenics\_optim} package will be an interesting perspective for future work. Providing interfaces to other interior-point solvers, especially open-source ones (CVXOPT, ECOS, Sedumi, etc.), is also planned for future releases.\\

The manuscript is organized as follows: section \ref{sec:conic-prog} introduces the conic programming framework and the concept of conic representable functions. The formulation and discretization of convex variational problems is discussed in section \ref{sec:var-prob} by means of a simple example. Section \ref{sec:Cheeger} discusses further aspects by considering a more advanced example. Finally, \ref{sec:gallery} provides a gallery of illustrative examples along with their formulation and some numerical results.

The \texttt{fenics\_optim} package can be downloaded from \url{https://gitlab.enpc.fr/navier-fenics/fenics-optim}. It contains test files as well as demo files corresponding to the examples discussed in the present paper.

\section{Conic programming framework}\label{sec:conic-prog}

\subsection{Conic programming in Mosek}\label{sec:cones}
The Mosek solver is dedicated to solving problems entering the \textit{conic programming framework} which can be written as:
\begin{equation}
\begin{array}{rl} \displaystyle{\min_\mathbf{x}} & \mathbf{c}\T \mathbf{x} \\
\text{s.t.} &\mathbf{b}_l \leq \mathbf{Ax}\leq \mathbf{b}_u\\
& \mathbf{x}\in \Kk
\end{array} \label{conic-programming}
\end{equation}
where vector $\mathbf{c}$ defines a linear objective functional, matrix $\mathbf{A}$ and vectors $\mathbf{b}_u,\mathbf{b}_l$ define linear inequality (or equality if $\mathbf{b}_u=\mathbf{b}_l$) constraints and where $\Kk=\Kk^1\times \Kk^2\times \ldots \times \Kk^p$ is a product of cones $\Kk^i \subset \RR^{d_i}$ so that $\mathbf{x}\in \Kk \Leftrightarrow \mathbf{x}^i \in \Kk^i \:\forall i=1,\ldots,p$ where $\mathbf{x}=(\mathbf{x}^1,\mathbf{x}^2, \ldots, \mathbf{x}^p)$. These cones can be of different kinds:
\begin{itemize}
\item $\Kk^i = \mathbb{R}^{d_i}$ i.e. no constraint on $\mathbf{x}^i$
\item $\Kk^i = (\mathbb{R}^+)^{d_i}$ is the positive orthant i.e. $\mathbf{x}^i \geq 0$
\item $\Kk^i = \Qq_{d_i}$ the quadratic Lorentz cone defined as:
\begin{equation}
\Qq_{d_i} = \{\mathbf{z} \in \RR^{d_i} \text{ s.t. } \mathbf{z}=(z_0, \bar{\mathbf{z}}) \text{ and } z_0 \geq \|\bar{\mathbf{z}}\|_2\}
\end{equation}
\item $\Kk^i = \Qq_{d_i}^r$ the rotated quadratic Lorentz cone defined as:
\begin{equation}
\Qq_{d_i}^r = \{\mathbf{z} \in \RR^{d_i} \text{ s.t. } \mathbf{z}=(z_0, z_1, \bar{\mathbf{z}}) \text{ and } 2z_0z_1 \geq \|\bar{\mathbf{z}}\|_2^2\}
\end{equation}
\item $\Kk^i = \Ss_{d_i}$ is the vectorial representation of the cone of semi-definite positive matrices $\mathbb{S}_{n_i}^+$ of dimension $n_i$ if $d_i=n_i(n_i+1)/2$ i.e.
\begin{align}
\Ss_{d_i}&= \{\vec(M) \text{ s.t. } M \in \mathbb{S}^+_{n_i}\} \\
\text{where } \mathbb{S}^+_{n_i}&= \{M \in \mathbb{R}^{n_i\times n_i} \text{ s.t. } M=M\T \text{ and } M\succeq 0 \} \notag
\end{align}
and in which the $\vec$ operator is the half-vectorization of a symmetric matrix obtained by collecting in a column vector the upper triangular part of the matrix. The elements are arranged in a column-wise manner and off-diagonal elements are multiplied by a $\sqrt{2}$ factor. For instance, for a $3\times 3$ matrix:
\begin{equation}
\vec(M) = (M_{1,1},\sqrt{2}M_{1,2},\sqrt{2}M_{1,3}, M_{2, 2}, \sqrt{2}M_{2,3}, M_{3,3})
\end{equation}
Note that this representation is such that $\vec(A)\T\vec(B) = \langle A, B\rangle_F$ where $\langle\cdot,\cdot\rangle_F$ denotes the Froebenius inner product of two matrices.
\end{itemize}

If $\Kk$ contains only cones of the first two kinds, then the resulting optimization problem \eqref{conic-programming} belongs to the class of Linear Programming (LP) problems. If, in addition, $\Kk$ contains quadratic cones $\Qq_{d_i}$ or $\Qq_{d_i}^r$, then the problem belongs to the class of Second-Order Cone Programming (SOCP) problems. Finally, when cones of the type $\Ss_{d_i}$ are present, the problem belongs to the class of Semi-Definite Programming (SDP) problems. Note that Quadratic Programming (QP) problems consisting of minimizing a quadratic functional under linear constraints can be seen as a particular instance of an SOCP problem as we will later discuss. 

There obviously exist dedicated algorithms for some classes of problem (e.g. the simplex method for LP, projected conjugate gradient methods for bound constrained QP, etc.). However, interior-point algorithms prove to be extremely efficient algorithms for all kinds of problems from LP up to difficult problems like SDP. It also turns out that a large variety of convex problems can be reformulated into an equivalent problem of the previously mentioned categories so that interior-point algorithms can be used to solve, in a robust and efficient manner, a large spectrum of convex optimization problems.

\subsection{Conic reformulations}\label{sec:conic-reformulations}
Most conic programming solvers other than Mosek (CVXOPT, Sedumi, SDPT3) use a default format similar to \eqref{conic-programming}. Aiming at optimizing a convex problem using a conic programming framework therefore requires a first reformulation step to fit into format \eqref{conic-programming}. In the following examples, we will consider a purely discrete setting in which optimization variables are in $\mathbb{R}^n$.

\subsubsection{$L^2$-norm constraint}
Let us consider the following $L^2$-norm constraint:
\begin{equation}
\|\mathbf{x}\|_2 \leq 1
\end{equation}
This can be readily observed to be the following quadratic cone constraint $(1,\mathbf{x})\in \Qq_{n+1}$. However, for this constraint to fit the general format of \eqref{conic-programming}, one must introduce an additional scalar variable $y$ such that the previous constraint can be equivalently written:
\begin{align}
y &= 1 \label{L2-ball}\\
\|\mathbf{x}\|_2 &\leq y \Leftrightarrow (y,\mathbf{x})\in\Qq_{n+1} \notag
\end{align}

\subsubsection{$L^1$-norm constraint}
Let us consider the following $L^1$-norm constraint:
\begin{equation}
\|\mathbf{x}\|_1 = \sum_{i=1}^n |x_i| \leq 1
\end{equation}
To reformulate this constraint, we introduce $n$ scalar auxiliary variables $y_i$ such that:
\begin{align}
\sum_{i=1}^n y_i &= 1 \\
|x_i| & \leq y_i \quad \forall i \notag
\end{align}
then each constraint with an absolute value can be written using two linear inequality constraints $x_i-y_i \leq 0$ and $0 \leq x_i+y_i$.

\subsubsection{Quadratic constraint}
Let us consider the case of a quadratic inequality constraint such as:
\begin{equation}
\frac{1}{2}\mathbf{x}\T \mathbf{Q}\mathbf{x} \leq b
\end{equation}

Matrix $\mathbf{Q}$ must necessarily be semi-definite positive for the constraint to be convex. In this case, introducing the Cholesky factor $\mathbf{C}$ of $\mathbf{Q}$ such that $\mathbf{Q}=\mathbf{C}\T\mathbf{C}$, one has:
\begin{equation}
\frac{1}{2}\mathbf{x}\T \mathbf{Q}\mathbf{x} = \frac{1}{2}\|\mathbf{C}\mathbf{x}\|_2^2 \leq b
\end{equation}
Introducing an auxiliary variable $\mathbf{y}$, the previous constraint can be equivalently reformulated as:
\begin{align}
\mathbf{C}\mathbf{x} - \mathbf{y} &= 0 \\
\|\mathbf{y}\|^2_2 &\leq 2b \notag
\end{align}
Finally adding two others scalar variables $z_0$ and $z_1$, we have:
\begin{align}
\mathbf{C}\mathbf{x} - \mathbf{y} &= 0 \\
z_0 &= b \notag \\
z_1 &= 1 \notag \\
\|\mathbf{y}\|^2_2 &\leq 2z_0z_1 \notag
\end{align}
where the last constraint is also the rotated quadratic cone constraint \mbox{$(z_0,z_1,\mathbf{y})\in\Qq_{n+2}^r$}

\subsubsection{Minimizing a $L^2$-norm}
Let us now consider the problem of minimizing the $L^2$-norm of $\mathbf{Bx}$ under some affine constraints:
\begin{equation}
\begin{array}{rl} \displaystyle{\min_\mathbf{x}} & \|\mathbf{Bx}\|_2 \\
\text{s.t.} & \mathbf{Ax} = \mathbf{b}\\
\end{array} \label{l2-norm}
\end{equation}
As such this problem does not fit \eqref{conic-programming} since the objective function is non-linear. In order to circumvent this, one needs to consider the epigraph of $F(\mathbf{x})=\|\mathbf{Bx}\|_2$ defined as $\epi F = \{(t,\mathbf{x}) \text{ s.t. } F(\mathbf{x}) \leq t\}$. Minimizing $F$ is then equivalent to minimizing $t$ under the constraint that $(t,\mathbf{x})\in \epi F$. For the present case, we therefore have:
\begin{equation}
\begin{array}{rl} \displaystyle{\min_{\mathbf{x},t}} & t \\
& \mathbf{Ax} = \mathbf{b}\\
& \|\mathbf{Bx}\|_2\leq t
\end{array}
\end{equation}
Introducing an additional variable $\mathbf{y}$ we have:
\begin{equation}
\begin{array}{rl} \displaystyle{\min_{\mathbf{x},\mathbf{y},t}} & t \\
& \mathbf{Ax} = \mathbf{b}\\
& \mathbf{Bx}-\mathbf{y} = 0\\
& \|\mathbf{y}\|_2\leq t
\end{array}
\end{equation}
where the last constraint is again a quadratic Lorentz cone constraint. Problem \eqref{l2-norm} is now a linear problem of the augmented optimization variables $(\mathbf{x},\mathbf{y},t)$ under linear and conic constraints.

\subsection{Conic representable sets and functions}

As previously mentioned, minimizing a convex function $F(\mathbf{x})$ can be turned into a linear problem with a convex non-linear constraint involving the epigraph of $F$. We will thus consider the class of \textit{conic representable functions} as the class of convex functions which can be expressed as follows:
\begin{equation}
\begin{array}{rl} \displaystyle{F(\mathbf{x}) = \min_{\mathbf{y}}} & \mathbf{c}_x\T\mathbf{x}+ \mathbf{c}_y\T\mathbf{y} \\
\text{s.t.} & \mathbf{b}_l\leq \mathbf{Ax}+\mathbf{By} \leq \mathbf{b}_u\\
& \mathbf{y}\in \Kk
\end{array} \label{conic-representable}
\end{equation}
in which $\Kk$ is again a product of cones of the kinds detailed in section \ref{sec:cones}.

For instance, consider the case of the $L^1$-norm, we have:
\begin{equation}
\begin{array}{rl} \displaystyle{\|\mathbf{x}\|_1 = \min_{\mathbf{y}\in\mathbb{R}^{n}}} &  \mathbf{e}\T\mathbf{y} \\
\text{s.t.} & \mathbf{0}\leq \mathbf{x}+\mathbf{y} \\
& \mathbf{x}-\mathbf{y} \leq \mathbf{0}
\end{array}
\end{equation}
where $\mathbf{e}=(1,\ldots, 1)$ whereas for the $L^2$-norm we have:
\begin{equation}
\begin{array}{rl} \displaystyle{\|\mathbf{x}\|_2 = \min_{\mathbf{y}\in\mathbb{R}^{n+1}}} &  y_0 \\
\text{s.t.} & \mathbf{x}-\mathbf{\bar{y}}=\mathbf{0} \\
& \mathbf{y} \in \Qq_{n+1}
\end{array} \label{L2norm-conic}
\end{equation}
where $\mathbf{y}=(y_0,y_1,\ldots,y_{n})$ and $\mathbf{\bar{y}} = (y_1,\ldots,y_{n})$. In this example, it can be seen that the representation \eqref{conic-representable} is not necessarily optimal in terms of number of additional variables, one could perfectly eliminate the $\mathbf{y}$ variable. However, in most practical cases, functions like the $L^2$-norm will quite often be composed with some linear operator (gradient, interpolation, etc.) so that introducing such additional variables will be necessary to fit format \eqref{conic-programming}.\\

Obviously, if $F$ is the indicator function of a convex set, then we have a similar notion of \textit{conic representable sets} for which only the constraints in \eqref{conic-representable} are relevant. 

\section{Variational problems and their discrete version}\label{sec:var-prob}
\subsection{A first illustrative example}\label{sec:obstacle}
Before describing the framework of variational formulation and discretization, let us first introduce a classical example of variational inequality, namely the obstacle problem. Let $\Omega$ be a bounded domain of $\RR^2$ and $u\in V = H_0^1(\Omega)$, $f\in H^{-1}(\Omega)$ and $g \in H^1(\Omega)\cap C^0$ such that $y\leq 0$ on $\partial \Omega$. The obstacle problem consists in solving:
\begin{equation}
\begin{array}{rl} \displaystyle{\inf_{u\in V}} & \displaystyle{\int_{\Omega}\dfrac{1}{2}\|\nabla u\|^2_2 \dx- \int_\Omega fu \dx} \\
\text{s.t.}  &u\in \Kk
\end{array} \label{obstacle}
\end{equation}
where $\Kk=\{v \in H_0^1(\Omega) \text{ s.t. } v \geq g \text{ on } \Omega\}$. Physically, this problem corresponds to that of a membrane described by an out-of-plane deflection $u$ and loaded by a vertical load $f$ which may potentially enter in contact with a rigid obstacle located on the surface $z=-g(x,y)$.

\subsection{Discretization}
Let us now consider some finite element discretization of $\Omega$ using a mesh $\Tt_h$ of $N_e$ triangular cells. For the displacement field $u$,  we consider a Lagrange piecewise linear interpolation represented by the discrete functional space $V_h = \{v \in C^0(\Omega) \text{ s.t. } v|_{T}\in \PP^1(T) \quad \forall T \in \Tt_h\}$ of dimension $N$. Interpolating the obstacle position $y$ on the same space $V_h$, a discrete approximation $\Kk_h$ of $\Kk$ consists in a pointwise inequality on the vectors $\mathbf{v},\mathbf{g}\in\RR^N$ of degrees of freedom of $v_h,g_h\in V_h$: $\Kk_h=\{v_h\in V_h \text{ s.t. } \mathbf{v} \geq \mathbf{g} \}$. Finally, introducing a quadrature formula with $M$ quadrature points for the first integral in \eqref{obstacle}, the discrete obstacle problem is now:
\begin{equation}
\begin{array}{rl} \displaystyle{\min_{\mathbf{u}\in \RR^N}} & \displaystyle{\sum_{g=1}^M \omega_g \dfrac{1}{2}\|\mathbf{B}_g \mathbf{u}\|^2_2 - \mathbf{f}\T\mathbf{u}} \\
\text{s.t.}  &\mathbf{u} \geq \mathbf{g}
\end{array} \label{obstacle-discr}
\end{equation}

In \eqref{obstacle-discr}, $\mathbf{B}_g\mathbf{u}\in \RR^2$  denotes the discrete gradient evaluated at the current quadrature point $g$ and $\omega_g$ is the associated quadrature weight. Note that since $u_h$ is linear, its gradient is piecewise-constant so that only one point per triangle $T$ with $\omega_g= |T|$ is sufficient for exact evaluation of the integral ($M=N_e$ in this case). Finally, $\mathbf{f}$ is the assembled finite-element vector corresponding to the linear form $L(u) = \int_\Omega fu \dx$.

The quadratic term in the objective function is now rewritten following section \ref{sec:conic-reformulations} as follows:
\begin{equation}
\begin{array}{rl} \displaystyle{\min_{\mathbf{u}\in \RR^N}} & \displaystyle{\sum_{g=1}^M \omega_g y_{g,0}  - \mathbf{f}\T\mathbf{u}} \\
\text{s.t.}  &\mathbf{u} \geq \mathbf{g} \\
& y_{g,1} = 1 \\
& \mathbf{B}_g\mathbf{u} - \begin{bmatrix}
y_{g,2} \\ y_{g,3}
\end{bmatrix} = 0 \qquad \forall g=1,\ldots, M\\
& (y_{g,0},y_{g,1},y_{g,2},y_{g,3}) \in \Qq_4^r
\end{array}
\end{equation}

Collecting the 4$M$ auxiliary variables $\mathbf{y}_g=(y_{g,0},y_{g,1},y_{g,2},y_{g,3})$ into a global vector $\mathbf{\widehat{y}}=(\mathbf{y}_1,\ldots\mathbf{y}_M) \in\RR^{4M}$, the previous problem can be rewritten as:
 \begin{equation}
\begin{array}{rl} \displaystyle{\min_{\mathbf{u}\in \RR^N, \mathbf{\widehat{y}}\in \RR^{4M}}} & \displaystyle{\mathbf{c}\T\mathbf{\widehat{y}} - \mathbf{f}\T\mathbf{u}} \\
\text{s.t.}  &\mathbf{u} \geq \mathbf{g} \\
& \mathbf{A}_u\mathbf{u} + \mathbf{A}_y\mathbf{\widehat{y}} = \mathbf{b}\\
& \mathbf{\widehat{y}} \in \Qq_4^r\times \cdots\times \Qq_4^r
\end{array}
\label{obstacle-discr-block}
\end{equation}
where $\mathbf{c} = (\omega_1, 0, 0, 0, \omega_2, 0 \ldots, \omega_M, 0, 0, 0)$, $\mathbf{A}_u = \begin{bmatrix}
0 \\ \mathbf{B}_1 \\ 0 \\ \mathbf{B}_2 \\ \vdots\\ 0 \\ \mathbf{B}_M
\end{bmatrix}$, $\mathbf{A}_y = \begin{bmatrix} -\mathbf{I} & & \\  & \ddots &  \\  &  & -\mathbf{I}\end{bmatrix}$ with $\mathbf{I} = \begin{bmatrix}
0 & 1 & 0 & 0 \\ 0 & 0 & 1 & 0 \\ 0 & 0 & 0 & 1
\end{bmatrix}$ and $\mathbf{b} = (1, 0, 0, 1, 0 \ldots, 1, 0, 0)$. This last formulation enables to see that problem \eqref{obstacle-discr} indeed fits into the general conic programming framework \eqref{conic-programming} but in a specific fashion since it possesses a block-wise structure induced by the quadrature rule. Indeed, each 4-dimensional block of auxiliary variables $y_g$ is decoupled from each other and is linked to the main unknown variable $\mathbf{u}$ through the evaluation of the discrete gradient at each point $g$. The conic reformulation performed in \eqref{obstacle-discr} is in fact the same for all quadrature points.

This observation motivates us to rewrite the initial continuous problem as:
\begin{equation}
\begin{array}{rl} \displaystyle{\inf_{u\in V}} & \displaystyle{\int_{\Omega}F(\nabla u) \dx- \int_\Omega fu \dx} \\
\text{s.t.}  &u\in \Kk
\end{array} \label{obstacle-cont-2}
\end{equation}
with $F(\mathbf{x}) = \dfrac{1}{2}\|\mathbf{x}\|^2_2$ which is conic representable as follows\footnote{Note that it would have been possible to work directly with function $\widetilde{F}(u):= \frac{1}{2}\|\nabla u\|_2^2$ which is also conic-representable}:
 \begin{equation}
\begin{array}{rl} \displaystyle{F(\mathbf{x}) = \min_{\mathbf{y}\in\mathbb{R}^{4}}} &  y_0 \\
\text{s.t.} & y_1 = 1\\
& \mathbf{x}- \begin{bmatrix}
y_{2} \\ y_{3}
\end{bmatrix}=\mathbf{0} \\
& \mathbf{y} \in \Qq_{4}^r
\end{array} \label{quad-representation}
\end{equation}

Introducing now the previously mentioned discretization and the quadrature formula, we aim at solving:
\begin{equation}
\begin{array}{rl} \displaystyle{\min_{\mathbf{u}\in \RR^N}} & \displaystyle{\sum_{g=1}^M \omega_g F(\mathbf{B}_g\mathbf{u})- \mathbf{f}\T\mathbf{u}} \\
\text{s.t.}  &\mathbf{u} \geq \mathbf{g}
\end{array} \label{obstacle-discr-2}
\end{equation}
which will be equivalent to \eqref{obstacle-discr} when injecting \eqref{quad-representation} into \eqref{obstacle-discr-2} since, for all $M$ evaluations of $F$, a 4-dimensional auxiliary vector $\mathbf{y}_g$ will be introduced as an additional minimization variable.

As a consequence, the \texttt{fenics\_optim} package has been particularly designed for a sub-class of problems of type \eqref{co-prob} in which $J(u) = \int_{\Omega} j(u)\dx$ for which integration will be handled by the FEniCS machinery and in which the user must specify the local conic representation of $j$.
 
\subsection{FEniCS formulation}\label{sec:fenics-obstacle}
In the following, we present the main part of a \texttt{fenics\_optim} script. More details on how an optimization problem is defined are discussed in \ref{appendix:structure}. In particular, it is possible to write \textit{manually} the discretized version of the obstacle problem based on \eqref{obstacle-discr-block} (see \ref{appendix:obstacle}). However, \texttt{fenics\_optim} also provides a more user-friendly way of modelling such problems which is based on \eqref{obstacle-cont-2} and \eqref{quad-representation} and will now be presented.

First, a simple unit square mesh and $\PP^1$ Lagrange function space \texttt{V} is defined using basic FEniCS commands. Homogeneous Dirichlet boundary conditions are also defined in variable \texttt{bc}. Finally, \texttt{obstacle} is the interpolant on $V$ of \[g(x,y) = g_0+a\sin(2\pi k_1 x)\cos(2\pi k_1 y)\sin(2\pi k_2 x)\cos(2\pi k_2 y)\]

In the following simulations, we took $g_0=-0.1$, $a=0.01$, $k_1=2$ and $k_2=8$. The loading is also assumed to be uniform and given by $f=-5$. The main part of the script starts by instantiating a \texttt{MosekProblem} object and adding a first optimization variable \texttt{u} living in the function space \texttt{V}, subject to Dirichlet boundary conditions \texttt{bc}. The \texttt{add\_var} method also enables to define a lower bound (resp. an upper bound) on an optimization variable by specifying a value for the \texttt{lx} (resp. \texttt{ux}) keyword. For the present case, we use \texttt{lx=obstacle} for enforcing $\mathbf{u} \geq \mathbf{g}$.
\begin{pythoncode}
prob = MosekProblem("Obstacle problem")
u = prob.add_var(V, bc=bc, lx=obstacle)

prob.add_obj_func(-dot(load,u)*dx)
\end{pythoncode}
where we also added the linear part of the objective function through the \texttt{add\_obj\_func} method.

The next step consists in defining the quadratic part of the objective function. For this purpose, we define a class inheriting from the base \texttt{MeshConvexFunction} class which must be instantiated by specifying on which previously defined optimization variable\footnote{see the notion of block-variables discussed in \ref{appendix:structure}.} this function will act (here the only possible variable is \texttt{u}). Moreover, we also specify the degree of the quadrature necessary for integrating the function (one-point quadrature used by default but written explicitly in the code snippet below). We must also define the \texttt{conic\_repr} method which will encode the conic representation \eqref{quad-representation}. We add a local optimization variable \texttt{Y} of dimension 4 which will belong to the cone \texttt{RQuad(4)} representing $\Qq_4^r$. Equality constraints are then added using the \texttt{add\_eq\_constraint} by specifying, as in \eqref{conic-representable}, a block matrix $\begin{bmatrix}\mathbf{A} & \mathbf{B}\end{bmatrix}$ and a right-hand side \texttt{b} ($=0$ by default). Note that both equality constraints could also have been written in a single one of row dimension 3. Finally, the local linear objective $(\mathbf{c}_x,\mathbf{c}_y)$ vector is defined using the \texttt{set\_linear\_term} method.
\begin{pythoncode}
class QuadraticTerm(MeshConvexFunction):
    def conic_repr(self, X):
        Y = self.add_var(4, cone=RQuad(4))
        self.add_eq_constraint([None, Y[1]], b=1)
        self.add_eq_constraint([X, -as_vector([Y[2], Y[3]])])
        self.set_linear_term([None, Y[0]])

F = QuadraticTerm(u, degree=0)
F.set_term(grad(u))
prob.add_convex_term(F)
\end{pythoncode}
Note that constraints and linear objectives are all defined in a block-wise manner, these blocks consisting of, first, the main variable which has been specified at instantiation (\texttt{u} in this case), then the additional local variables (\texttt{Y} here). Besides, these blocks are represented in terms of their action on the block variables using UFL expressions. 

The \texttt{set\_term} method enables to evaluate $F$ for the gradient of $u$ using the UFL \texttt{grad} operator. This function is then added to the global optimization problem. Finally, optimization (minimization by default) is performed by calling the \texttt{optimize} method of the \texttt{MosekProblem} object:
\begin{pythoncode}
prob.optimize()
\end{pythoncode}

\begin{figure}
\begin{center}
\begin{subfigure}[b]{0.57\textwidth}
\includegraphics[width=\textwidth]{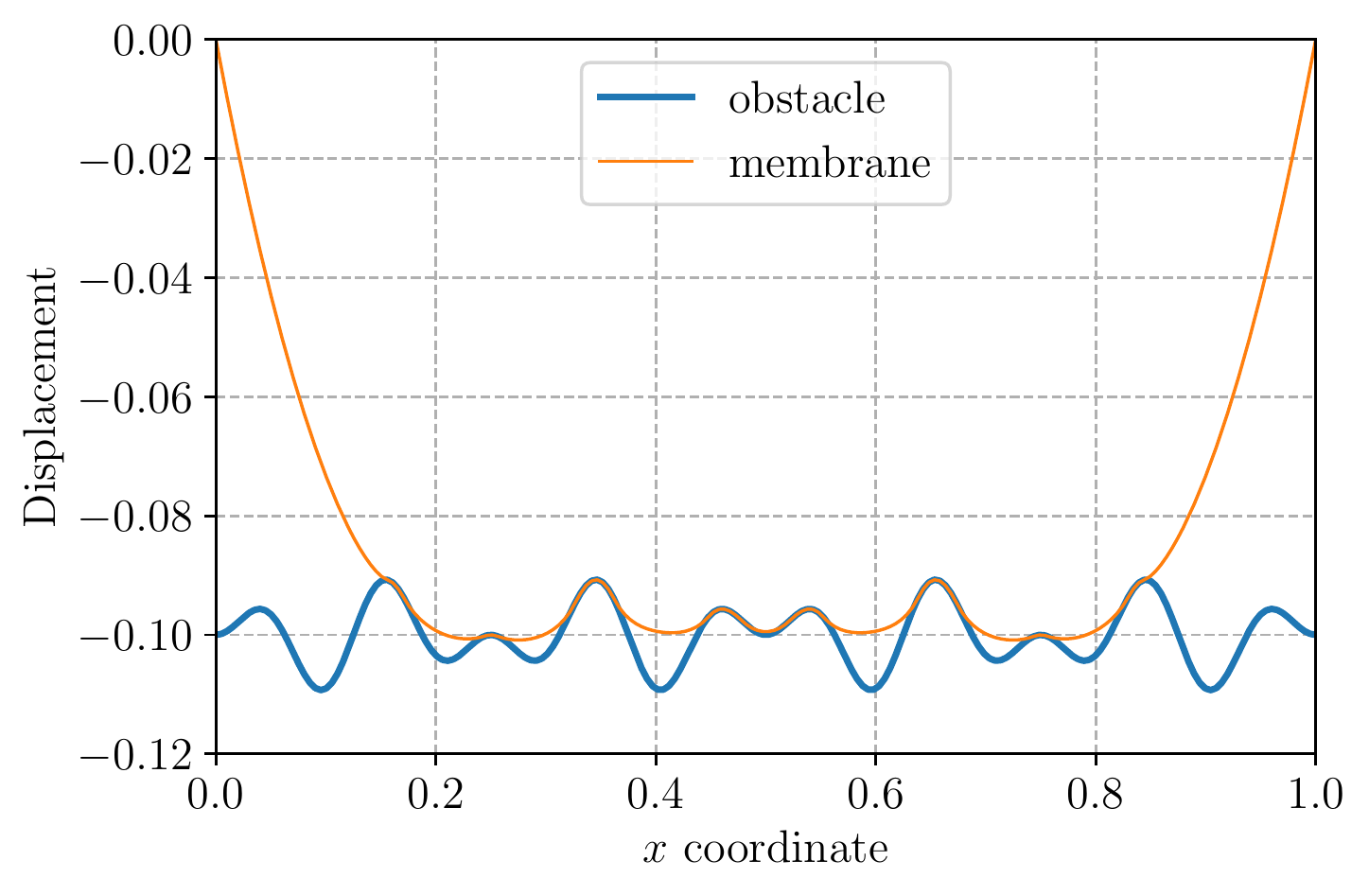}
\caption{Membrane displacement in contact with the obstacle along $y=0.5$}
\end{subfigure}
\hfill
\begin{subfigure}[b]{0.42\textwidth}
\includegraphics[width=\textwidth]{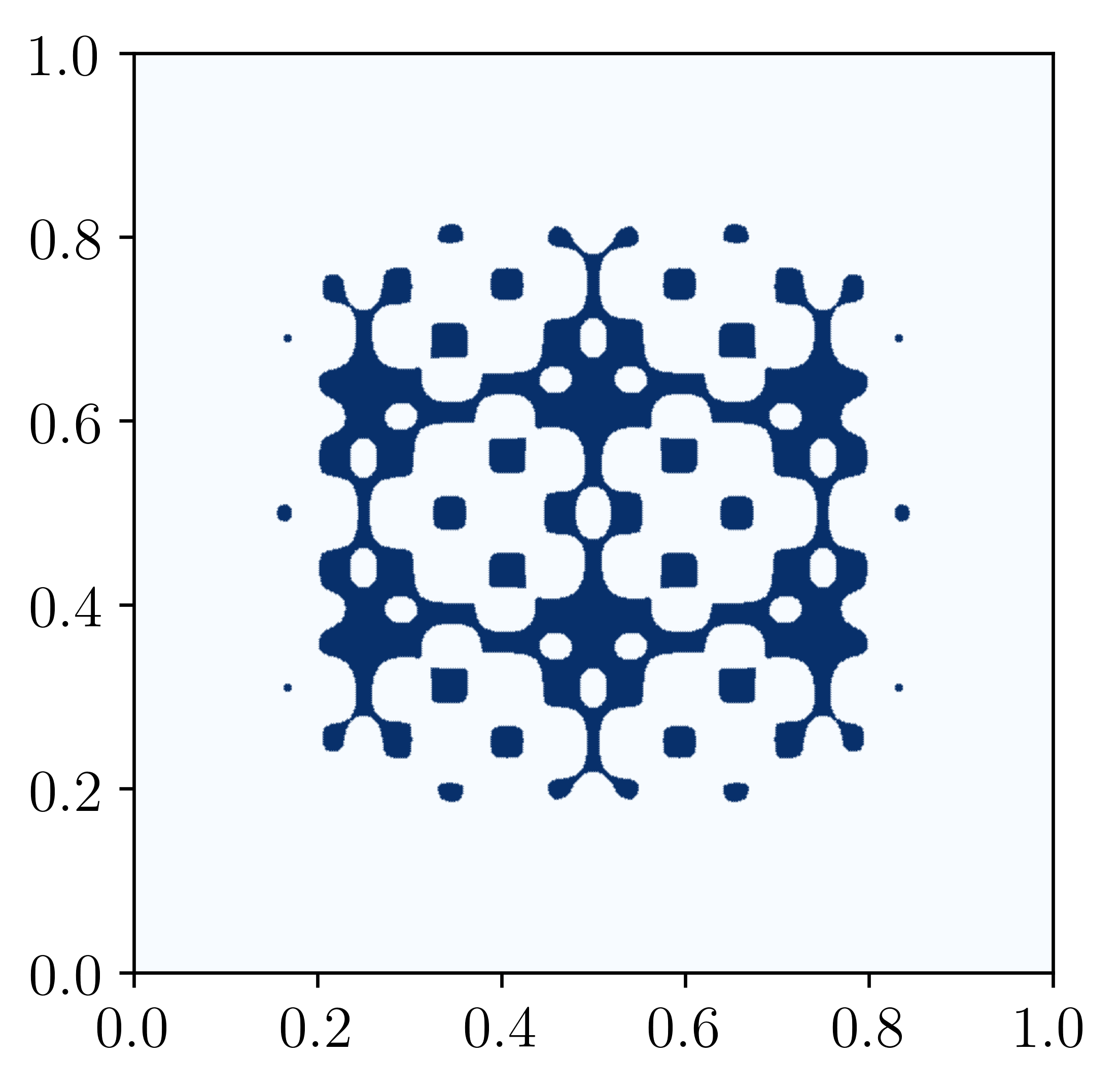}
\caption{Contact area in blue}
\end{subfigure}
\end{center}
\caption{Results of the obstacle problem for $h=1/400$}
\label{obstacle-fig}
\end{figure}

For validation and performance comparison, the obstacle problem has been solved for various mesh sizes using the \texttt{fenics\_optim} toolbox as well as using PETSc's TAO quadratic bound-constrained solver \cite{petsc-web-page,munson2012tao} which is particularly well suited for this kind of problems. We used the Trust Region Newton Method (TRON) and an ILU-preconditioned conjugate gradient solver for the inner iterations. Results in terms of optimal objective function value, total optimization time and number of iterations have been reported for both methods in table \ref{tab:obstacle-results}. Note that default convergence tolerances have been used in both cases and that total optimization time includes the presolve step of Mosek which can efficiently eliminate redundant linear constraints for instance. It can be observed that both approach yield close results in terms of optimal objective values and that TAO's solver is more efficient than Mosek in terms of optimization time as expected, mainly because of the small number of iterations needed to reach convergence but also because no additional variables are introduced when using TAO. However, Mosek surprisingly becomes quite competitive for large-scale problems because of its number of iterations scaling quite weakly with the problem size, contrary to the TRON algorithm. 
Membrane displacement along the line $y=0.5$ and contact area for $h=1/400$ have been represented in Figure \ref{obstacle-fig}.

\begin{table}
\begin{tabular}{|c||c|c|c||c|c|c|}
\hline 
\multirow{2}{*}{Mesh size} & \multicolumn{3}{c||}{\textbf{Interior point (Mosek)}} & \multicolumn{3}{c|}{\textbf{TRON algorithm (TAO)}} \\
\cline{2-7} 
 & Objective & Opt. time & iter. & Objective & Opt. time & iter.\\ 
\hline 
$h=1/25$ & -0.265081 & \phantom{4}0.13 s & 14 & -0.265082 & \phantom{4}0.09 s & 5 \\ 
\hline 
$h=1/50$ & -0.264932 & \phantom{4}0.56 s & 15 & -0.264932 & \phantom{4}0.22 s & 6 \\ 
\hline 
$h=1/100$ & -0.264883 & \phantom{4}2.27 s & 16 &  -0.264884 & \phantom{4}1.04 s & 10\\ 
\hline 
$h=1/200$ & -0.264867 & 10.04 s & 19 & -0.264871 & \phantom{4}6.03 s & 14 \\ 
\hline 
$h=1/400$ & -0.264864 & 48.95 s & 20 & -0.264868 &  47.79 s & 22  \\ 
\hline 
\end{tabular} 
\caption{Comparison between the \texttt{fenics\_optim} implementation of the obstacle problem relying on the interior point Mosek solver and TAO's bound-constrained TRON solver.}
\label{tab:obstacle-results}
\end{table}

\section{A more advanced example}\label{sec:Cheeger}
Let us now consider the following problem:
\begin{equation}
\begin{array}{rl} \displaystyle{\inf_{u\in V}} & \displaystyle{\int_{\Omega}\|\nabla u\|_2 \dx} \\
\text{s.t.}  & \int_\Omega fu \dx = 1
\end{array} \label{Cheeger}
\end{equation}
This problem is known to be related to antiplane limit analysis problems in mechanics as well as to the Cheeger problem and the eigenvalue fo the 1-Laplacian when $f=1$ \cite{cheeger1969lower,carlier2009approximation,carlier2011projection}. In this particular case, the solution of \eqref{Cheeger} can indeed be shown to be proportional to the characteristic function of a subset $C_\Omega\subseteq \Omega$ known as the Cheeger set of $\Omega$ which is the solution of:
\begin{equation}
C_{\Omega}:=\begin{array}{rl} \displaystyle{\argmin_{\omega \subseteq \Omega}} & \displaystyle{\dfrac{|\partial \omega|}{|\omega|}}
\end{array} 
\end{equation}
that is the subset minimizing the ratio of perimeter over area, the associated optimal value of this ratio $c_\Omega$ being known as the \textit{Cheeger constant}.

This problem is not strictly convex and is particularly difficult to solve using standard algorithms due to the highly non-smooth objective term. Again, introducing a $\PP^k$ Lagrange discretization for $u$, we aim at solving the following discrete problem:
\begin{equation}
\begin{array}{rl} \displaystyle{\min_{\mathbf{u}\in\RR^N}} & \displaystyle{\sum_{g=1}^M\omega_g F(\mathbf{B}_g \mathbf{u})} \\
\text{s.t.}  & \mathbf{f}\T\mathbf{u} = 1
\end{array} \label{Cheeger-discr}
\end{equation}
where $F(\mathbf{x})=\|\mathbf{x}\|_2$ with its conic representation being given by \eqref{L2norm-conic}. Similarly to the obstacle problem, choosing a $\PP^1$ discretization requires only a one-Gauss point quadrature rule for the objective function evaluation. For $\PP^k$ with $k\geq 2$, the quadrature is always inexact and Gaussian quadrature is not necessarily optimal. For the particular case $k=2$, one can choose a \textit{vertex} quadrature scheme on the simplex triangle to ensure that the discrete integral is approximated by excess:
\begin{equation}
\int_T \|\mathbf{r}(x,y)\|\dx \lesssim \dfrac{|T|}{3}\sum_{i=1}^3  \|\mathbf{r}(x_i,y_i)\| \label{vertex}
\end{equation}
where $(x_i,y_i)$ denote the simplex vertices. The choice of the quadrature scheme can also be made when defining the corresponding \texttt{MeshConvexFunction}:
\begin{pythoncode}
class L2Norm(MeshConvexFunction):
    """ Defines the L2-norm function ||x||_2 """
    def conic_repr(self, X):
        d = self.dim_x
        Y = self.add_var(d+1, cone=Quad(d+1))
        Ybar = as_vector([Y[i] for i in range(1, d+1)])
        self.add_eq_constraint([X, -Ybar])
        self.set_linear_term([None, Y[0]])
        
prob = MosekProblem("Cheeger problem")
u = prob.add_var(V, bc=bc)
    
if degree == 1:
    F = L2Norm(u)
elif degree == 2:
	F = L2Norm(u, "vertex")
else:
	F = L2Norm(u, degree = degree)
F.set_term(grad(u))
prob.add_convex_term(F)
\end{pythoncode}

In the previous code, \texttt{degree} denotes the polynomial degree $k$ of function space \texttt{V}. If $k=1$, the default one-point quadrature rule is used, if $k=2$ the above-mentioned \textit{vertex} scheme is used, otherwise a default Gaussian quadrature rule for polynomials of degree $k$ is used. \texttt{Quad(d+1)} corresponds to the quadratic Lorentz cone $\Qq_{d+1}$ of dimension $d+1$ where $d=\texttt{self.dim\_x}$ is the dimension of the \texttt{X} variable.

In the Cheeger problem, a normalization constraint must also be added. This can again be done by adding a convex term including only the corresponding constraint or it can also be added directly to the \texttt{MosekProblem} instance by defining the function space for the Lagrange multiplier corresponding to the constraint (here it is scalar so we use a \texttt{"Real"} function space) and passing the corresponding constraint in its weak form as follows:
\begin{pythoncode}
f = Constant(1.)
R = FunctionSpace(mesh, "Real", 0)
def constraint(l):
    return [l*f*u*dx]
prob.add_eq_constraint(R, A=constraint, b=1)
\end{pythoncode}

\subsection{Discontinuous Galerkin discretization}
Problem \eqref{Cheeger} can be discretized using standard Lagrange finite elements but also using Discontinous Galerkin discretization, in this case the gradient $L^2$-norm objective term is completed by absolute values of the jumps of $u$:
\begin{equation}
\begin{array}{rl} \displaystyle{\min_{u\in V}} & \displaystyle{\int_{\Omega}\|\nabla u\|_2 \dx + \int_{\Gamma} |\jump{u}|\dS + \int_{\Gamma_D} |u|\dS} \\
\text{s.t.}  & \int_\Omega fu \dx = 1
\end{array} \label{Cheeger-disc}
\end{equation}
where $\Gamma$ denotes the set of internal edges, $\Gamma_D$ the Dirichlet boundary part and $\jump{u}=u^+-u^-$ is the jump across $\Gamma$.

The discretized version using discontinuous $\PP^k_d$ Lagrange finite elements reads as:
\begin{equation}
\begin{array}{rl} \displaystyle{\min_{\mathbf{u}\in\RR^N}} & \displaystyle{\sum_{g=1}^M\omega_g F(\mathbf{B}_g \mathbf{u}) + \sum_{g_e=1}^{M_e} \omega_{g_e}G(\mathbf{J}_{g_e}\mathbf{u}) + \sum_{g_d=1}^{M_d} \omega_{g_d}G(\mathbf{T}_{g_d}\mathbf{u})} \\
\text{s.t.}  & \mathbf{f}\T\mathbf{u} = 1
\end{array} \label{Cheeger-discr-DG}
\end{equation}
where $G(x) = |x|$, $g_e$ (resp. $g_d$) denotes a current quadrature point on the internal (resp. Dirichlet) facets, $M_e$ (resp. $M_d$) denoting the total number of such points and $\omega_{g_e}$ (resp. $\omega_{g_d}$) the associated quadrature weights. Finally, $\mathbf{J}_{g_e}\mathbf{u}$ denotes the evaluation of $\jump{u}$ at the quadrature point $g_e$ and $\mathbf{T}_{g_d}\mathbf{u}$ the evaluation of $u$ at $g_d$.

Similarly to the previously introduced \texttt{MeshConvexFunction}, we define a \texttt{FacetConvexFunction} corresponding to the conic representable convex function $G(x)$:
\begin{pythoncode}
class AbsValue(FacetConvexFunction):
    def conic_repr(self, X):
        Y = self.add_var()
        self.add_ineq_constraint(A=[X, -Y], bu=0)
        self.add_ineq_constraint(A=[-X, -Y], bu=0)
        self.set_linear_term([None, Y])
\end{pythoncode}

When instantiating such a \texttt{FacetConvexFunction}, integration will (by default) be performed both on internal edges (in FEniCS the corresponding integration measure symbol is \texttt{dS}) and on external edges (FEniCS symbol being \texttt{ds}). If the Dirichlet boundary does not cover the entire boundary, then the \texttt{ds} measure can be restricted to the corresponding part. Again, the optimization variable on which acts the function must be specified and the desired quadrature rule can also be passed as an argument when instantiating the function. The \texttt{set\_term} method can now take a list of UFL expression associated to the different integration measures. In the present case, $G$ is evaluated for $\jump{u}$ on \texttt{dS} and for $u$ on \texttt{ds}:
\begin{pythoncode}
G = AbsValue(u)
G.set_term([jump(u), u])
prob.add_convex_term(G)
\end{pythoncode}
By default, facet integrals are evaluated using the \textit{vertex} scheme.

\subsection{Numerical example}
\begin{figure}
\begin{center}
\begin{subfigure}{0.65\textwidth}
\includegraphics[width=\textwidth]{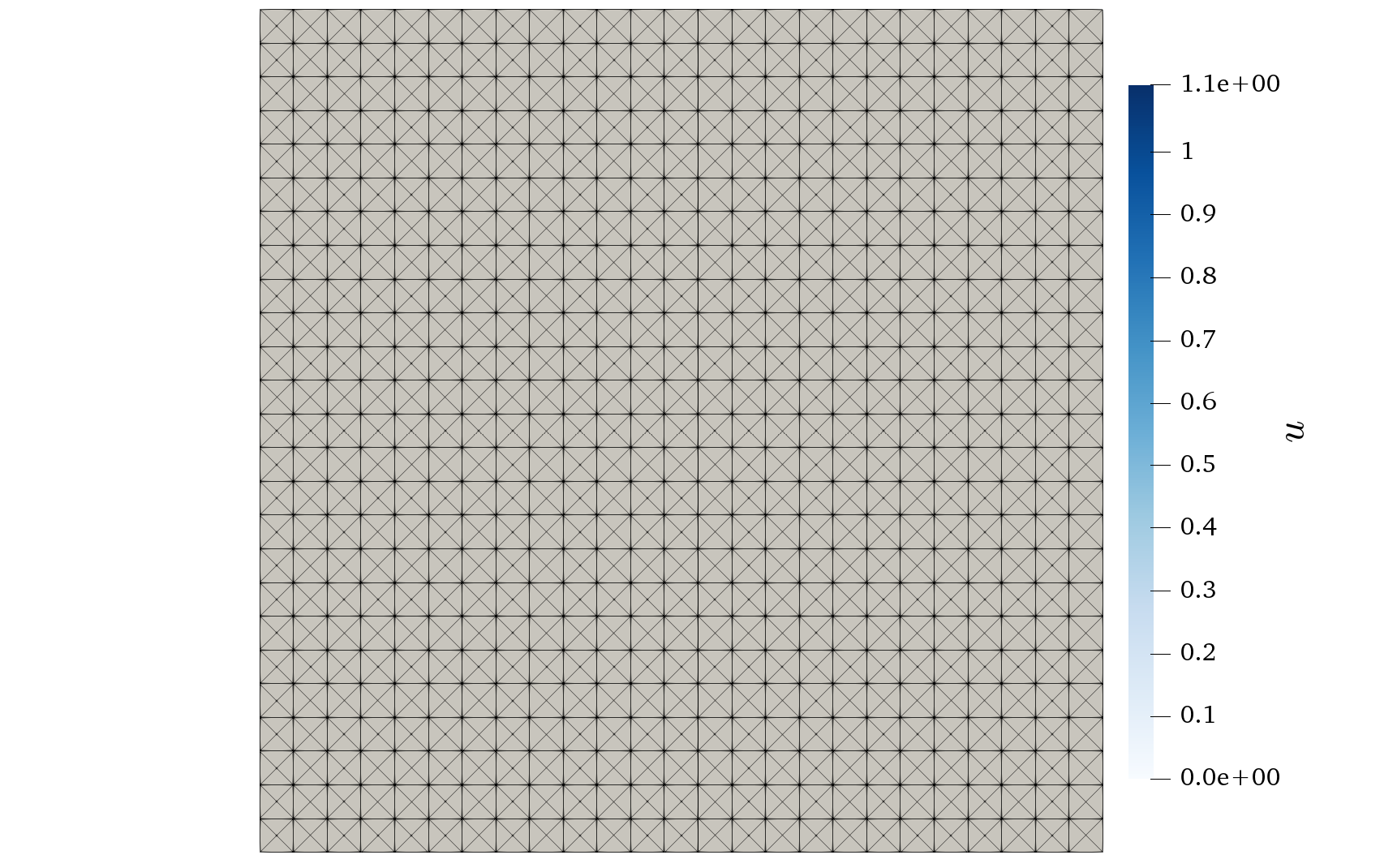}
\caption{25$\times$25 mesh}
\end{subfigure}
\\
\begin{subfigure}{0.45\textwidth}
\includegraphics[width=\textwidth]{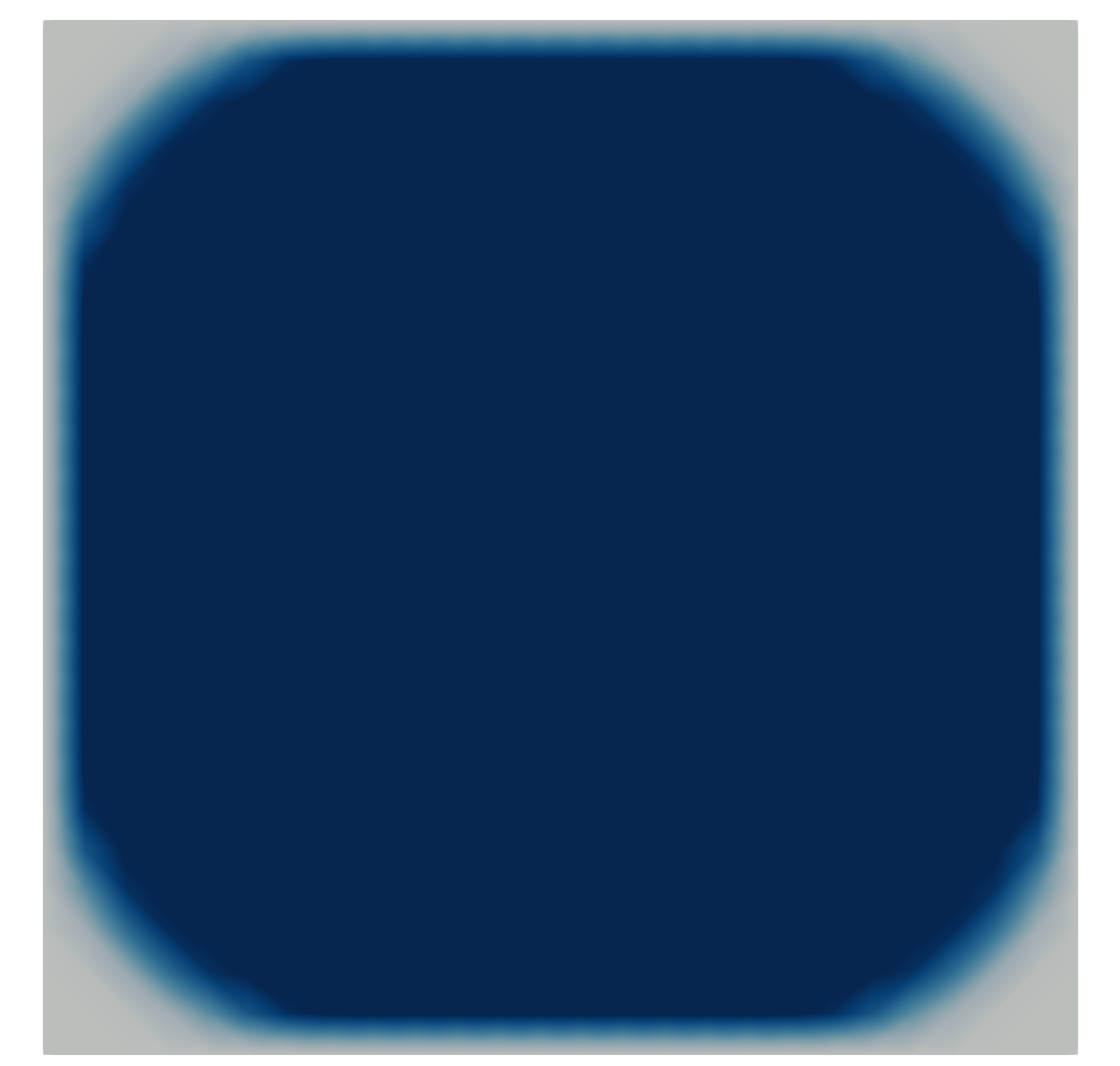}
\caption{CG-1}
\end{subfigure}
\hfill
\begin{subfigure}{0.45\textwidth}
\includegraphics[width=\textwidth]{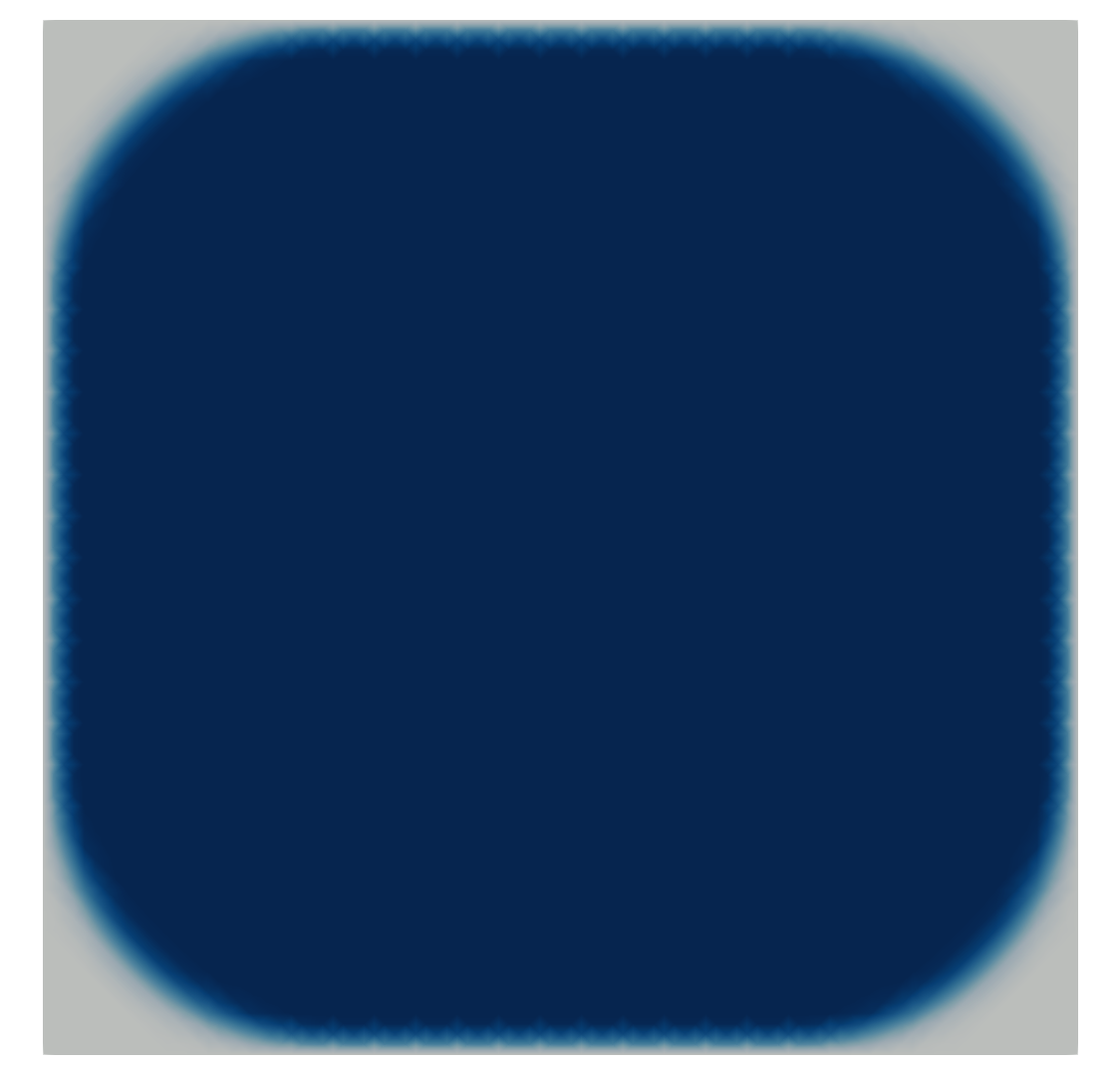}
\caption{CG-2}
\end{subfigure}\\
\begin{subfigure}{0.45\textwidth}
\includegraphics[width=\textwidth]{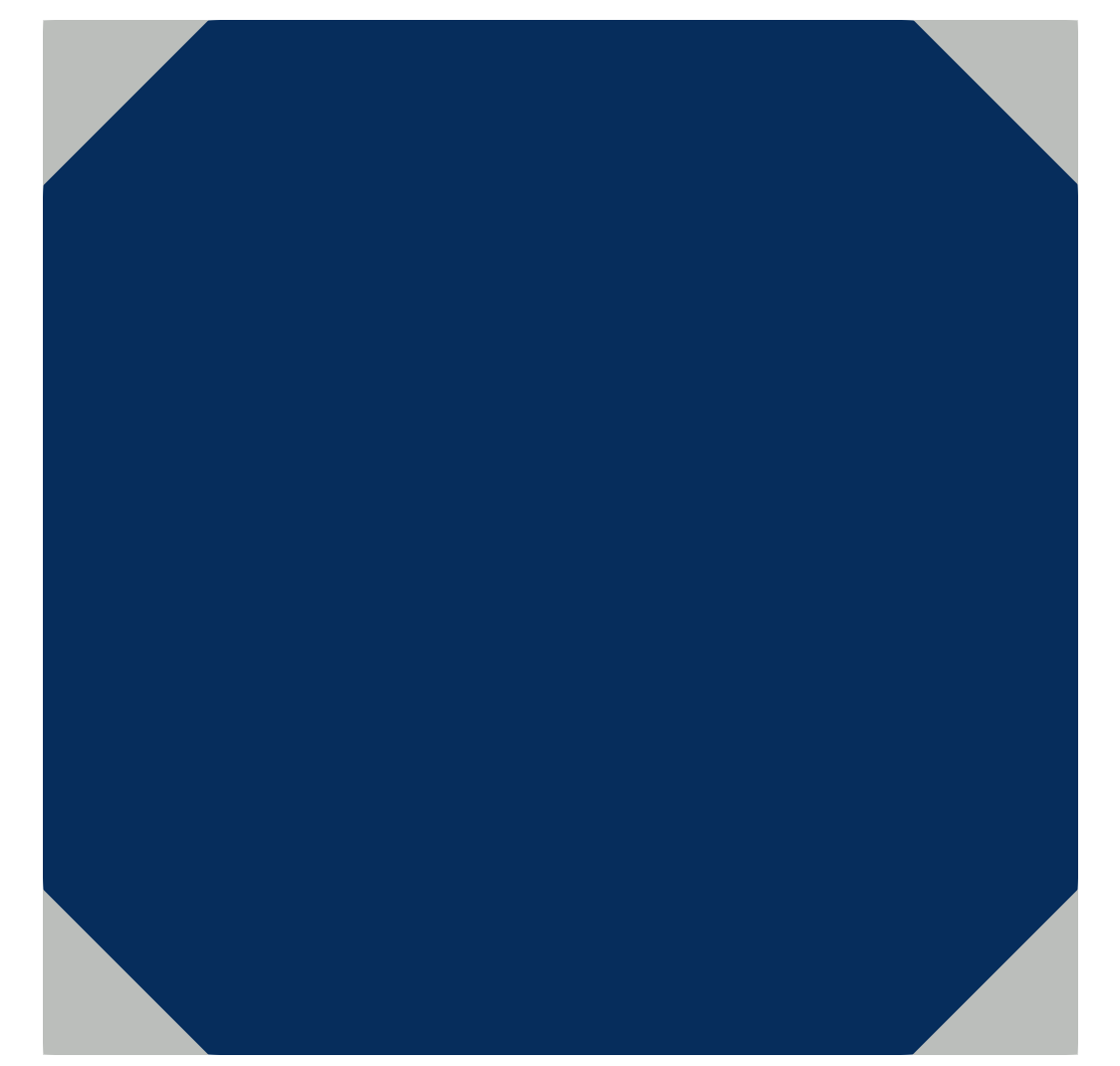}
\caption{DG-0}
\end{subfigure}
\hfill
\begin{subfigure}{0.45\textwidth}
\includegraphics[width=\textwidth]{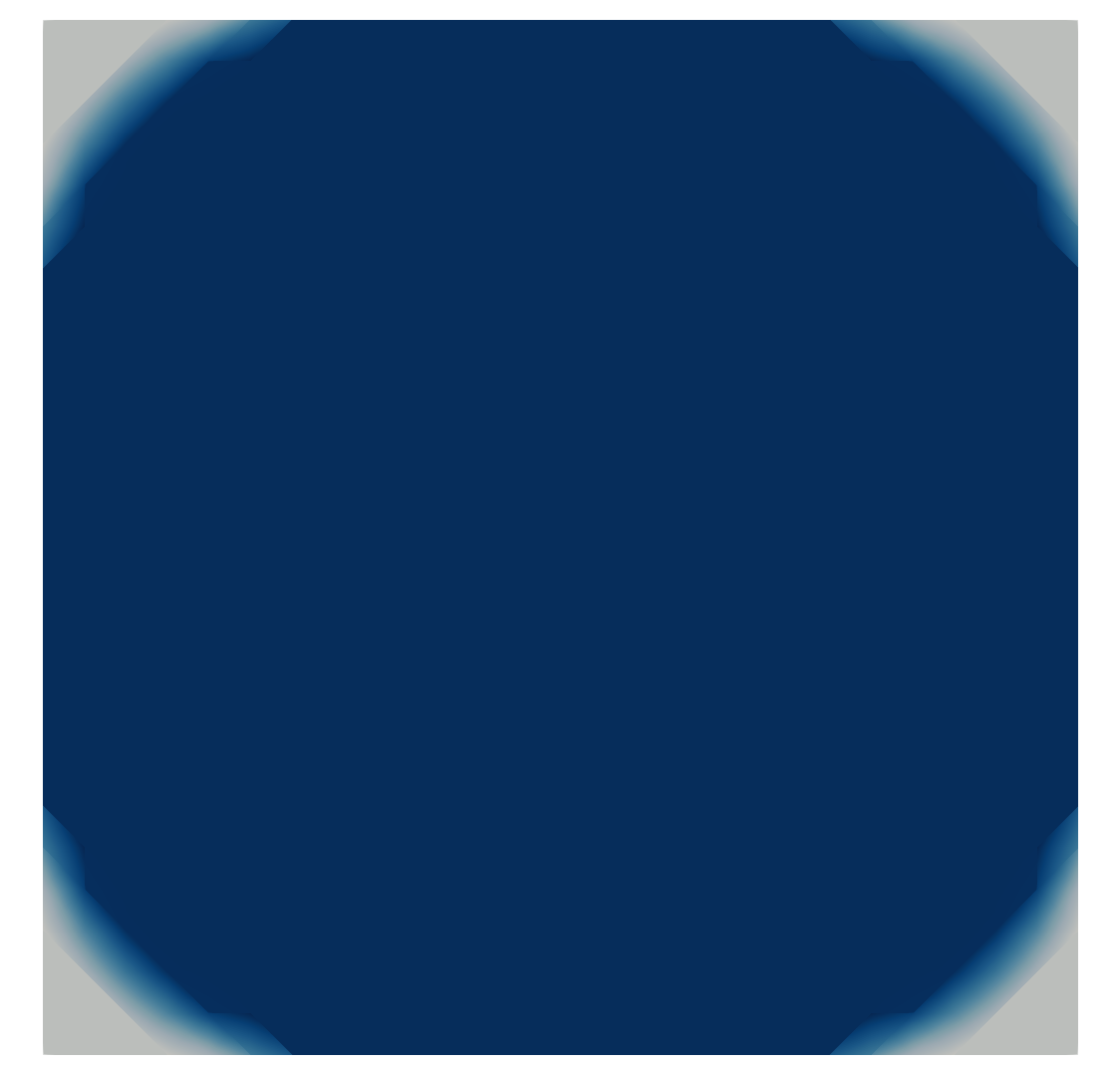}
\caption{DG-1}
\end{subfigure}
\end{center}
\caption{Results of the Cheeger problem for various discretizations on the unit square: continuous Galerkin (CG) and discontinuous Galerkin (DG) of degrees $k=0,1$ or $2$.}
\label{Cheeger-fields}
\end{figure}

We consider the problem of finding the Cheeger set of the unit square $\Omega=[0;1]^2$. The exact solution of this problem is known to be the unit square rounded by circles of radius $\rho = \dfrac{1}{2+\sqrt\pi}$ in its four corners, the associated Cheeger constant being $c_\Omega = 1/\rho$ \cite{strang1979minimax,overton1985numerical}. Results of the optimal field $u$ for various discretization schemes have been represented on Figure \ref{Cheeger-fields}. For all the retained discretization choices, the obtained Cheeger constant estimates are necessarily upper bounds to the exact one, in particular because of the choice of \textit{vertex} quadrature schemes ensuring upper bound estimations such as \eqref{vertex}. It can be seen on Figure \ref{Cheeger-fields} that all schemes yield a correct approximation of the Cheeger set, except for the DG-0 scheme which is too stiff and produces straight edges in the corners, following the structured mesh edges.

\subsection{A $H(\div)$-conforming discretization for the dual problem}
It can be easily shown through Fenchel-Rockafellar duality that problem \eqref{Cheeger} is equivalent to the following dual problem (see \cite{carlier2009approximation} for instance):
\begin{equation}
\begin{array}{rl} \displaystyle{\sup_{\lambda\in \RR, \bsig\in W}} & \displaystyle{\lambda} \\
\text{s.t.}  & \lambda f = \div\bsig \quad \text{in }\Omega \\
& \|\bsig\|_2 \leq 1
\end{array} \label{Cheeger-dual}
\end{equation}

A natural discretization strategy for such a problem is to use $H(\div)$-conforming elements such as the Raviart-Thomas element. Here, we will use the lowest Raviart-Thomas element, noted $RT_1$ by the FEniCS definition \cite{logg2012automated}. For the \texttt{fenics\_optim} implementation, two minimization variables are defined: $\lambda$ belonging to a scalar \texttt{"Real"} function space and $\bsig \in RT_1$.
Since for $\bsig \in RT_1$, $\div \bsig \in \PP^0$, we write the constraint equation using $\PP^0$ Lagrange multipliers:
\begin{pythoncode}
N = 50
mesh = UnitSquareMesh(N, N, "crossed")

VRT = FunctionSpace(mesh, "RT", 1)
R = FunctionSpace(mesh, "Real", 0)
VDG0 = FunctionSpace(mesh, "DG", 0)

prob = MosekProblem("Cheeger dual")
lamb, sig = prob.add_var([R, VRT])

f = Constant(1.)
def constraint(u):
    return [lamb*f*u*dx, -u*div(sig)*dx]
prob.add_eq_constraint(VDG0, A=constraint, name="u")
\end{pythoncode}

Finally, since $\bsig \in \PP^1$ on a triangle, if the constraint $\|\bsig\|_2$ is satisfied at the three vertices, it is satisfied everywhere by convexity. We here define a \texttt{MeshConvexFunction} representing the characteristic function of a $L^2$-ball constraint and select the \texttt{"vertex"} quadrature scheme so that the constraint will be indeed satisfied at the three vertices. Finally, the objective function is defined through the \texttt{add\_obj\_func} method of the problem instance:
\begin{pythoncode}
class L2Ball(MeshConvexFunction):
    """ Defines the L2-ball constraint ||x||_2 <= 1 """
    def conic_repr(self, X):
        d = self.dim_x
        Y = self.add_var(d+1, cone=Quad(d+1))
        Ybar = as_vector([Y[i] for i in range(1, d+1)])
        self.add_eq_constraint([X, -Ybar])
        self.add_eq_constraint([None, Y[0]], b=1)
        
F = L2Ball(sig, "vertex")
F.set_term(sig)
prob.add_convex_term(F)

prob.add_obj_func([1, None])
\end{pythoncode}

With the above-mentioned discretization and quadrature choice, it can easily be shown that the discrete version of \eqref{Cheeger-dual} will produce  a lower bound of the exact Cheeger constant. For instance, for a $25\times 25$ mesh, we obtained:
\begin{equation}
c_\Omega^{RT_1} \approx 3.704 \leq c_\Omega \approx 3.772 \leq c_\Omega^{DG1} \approx 3.800
\end{equation}

\begin{figure}
\begin{center}
\includegraphics[width=0.6\textwidth]{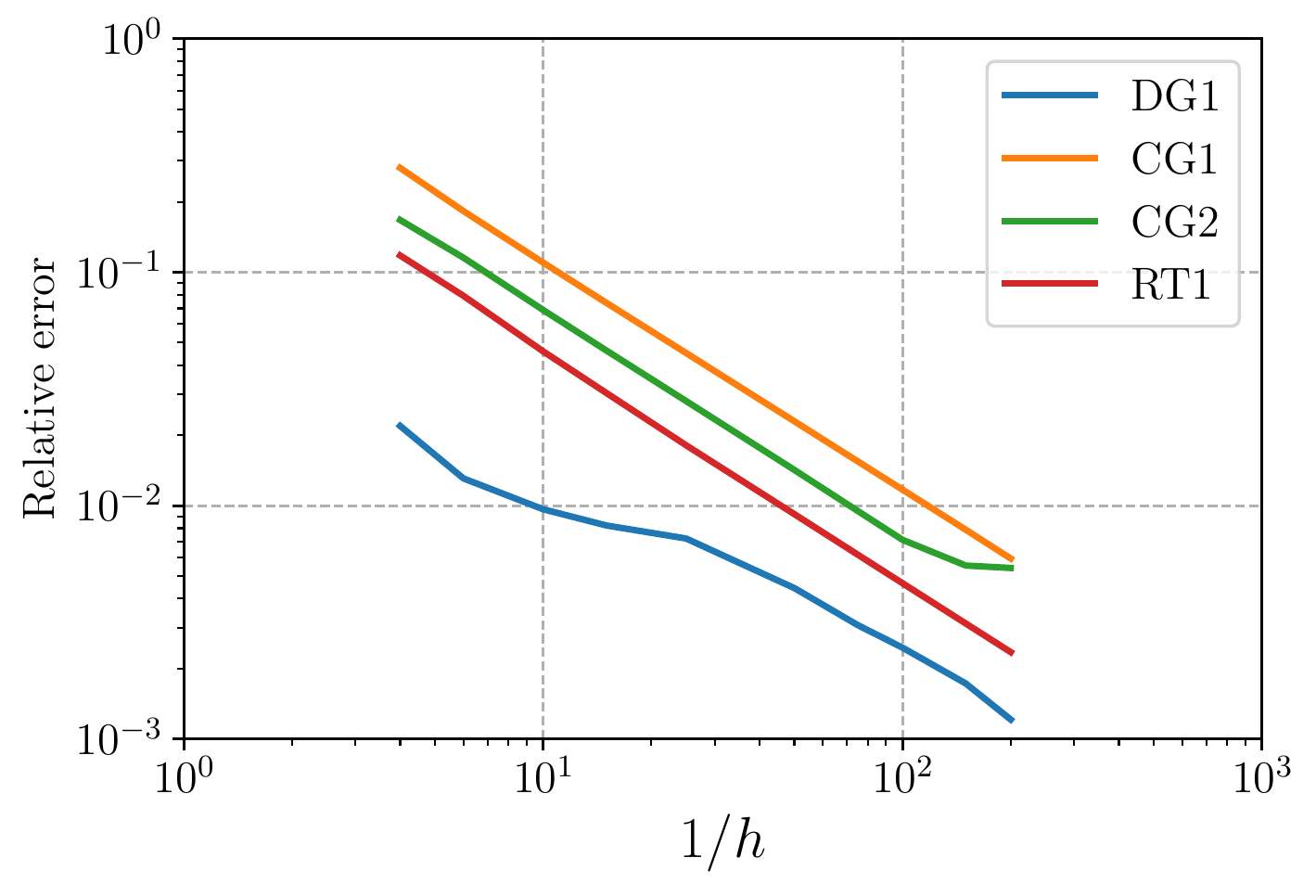}
\end{center}
\caption{Convergence results on the Cheeger problem}
\label{cheeger-convergence}
\end{figure}

\begin{figure}
\begin{center}
\includegraphics[width=0.4\textwidth]{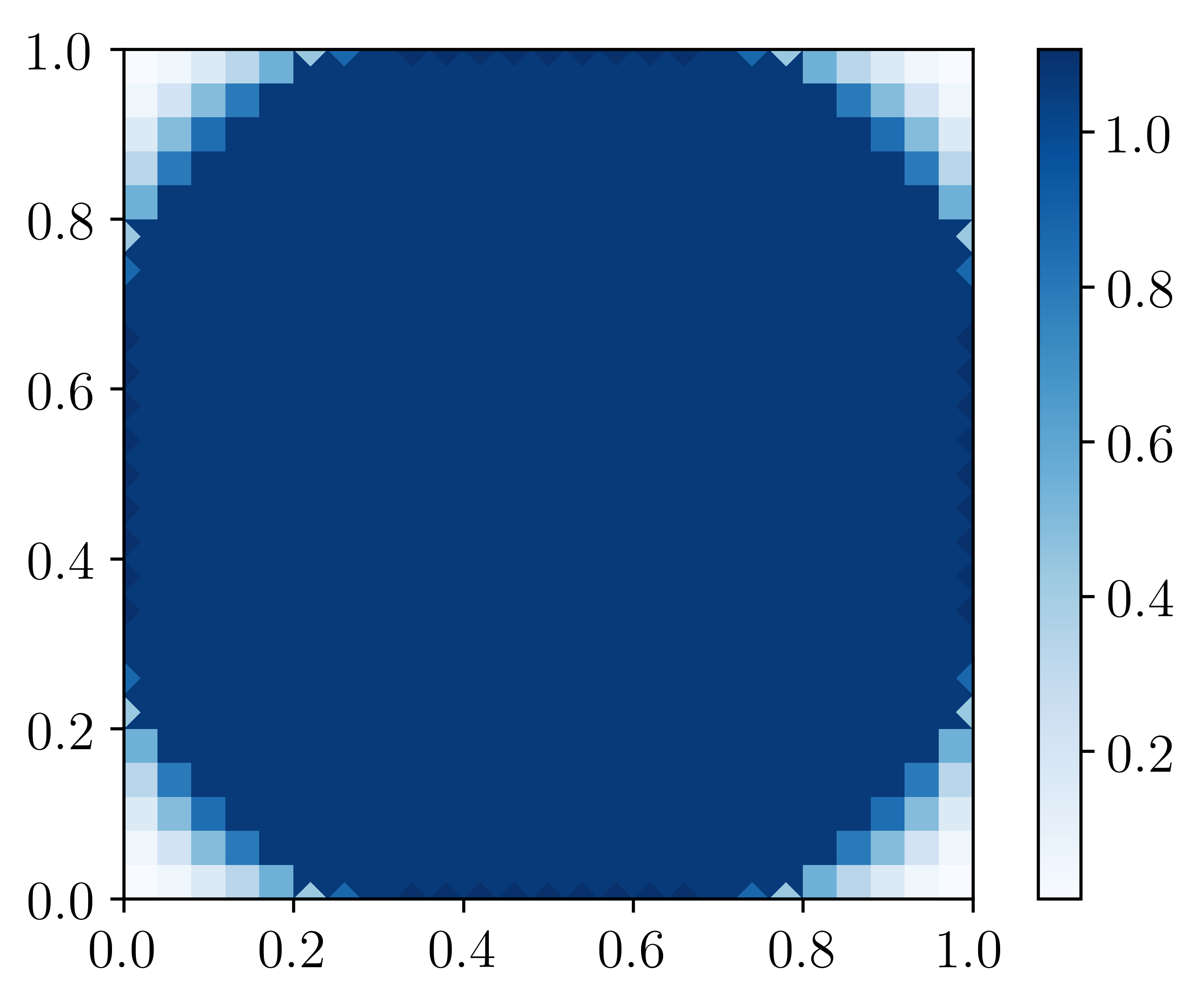}
\end{center}
\caption{Optimal $u$ field from the RT discretization}
\label{cheeger-dual-field}
\end{figure}

Convergence results of the numerical Cheeger constant estimate $c_{\Omega,h}$ obtained with the previous CG/DG discretizations as well as with the present RT discretization have been reported in Figure \ref{cheeger-convergence}. The relative error is computed as $\epsilon (c_{\Omega,h}/c_\Omega-1)$ where $\epsilon=-1$ for the RT discretization and $\epsilon=1$ otherwise. We observe in particular that the DG1 scheme is the most accurate and that all schemes have the same convergence rate in $O(h)$. Finally, primal-dual solvers such as Mosek also provide access to the optimal values of constraint Lagrange multipliers. The Lagrange multiplier associated with the constraint $\lambda f = \div\bsig$ can be interpreted as the field $u$ from the primal problem. This Lagrange multiplier, which belongs to a DG0 space, has been represented in Figure \ref{cheeger-dual-field}.


\subsection{A library of convex representable functions}\label{sec:libary}
In the \texttt{fenics\_optim} library, instead of defining each time the conic representation of usual functions, a library of common convex functions has been already implemented, including:
\begin{itemize}
\item linear functions $F(\mathbf{x}) = \mathbf{c}\T\mathbf{x}$
\item quadratic functions $F(\mathbf{x}) = \frac{1}{2}(\mathbf{x}-\mathbf{x}_0)\T\mathbf{Q}\T\mathbf{Q}(\mathbf{x}-\mathbf{x}_0)$
\item absolute value $F(x)=|x|$
\item $L^1$, $L^2$ and $L^\infty$ norms
\item $L^1$, $L^2$ and $L^\infty$ balls characteristic functions
\end{itemize}

\begin{figure}
\begin{center}
\begin{subfigure}{0.32\textwidth}
\includegraphics[width=\textwidth]{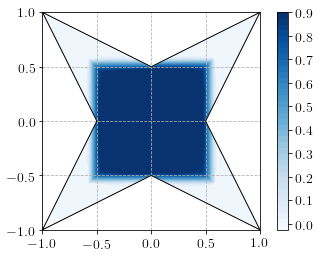}
\caption{$L^1$-norm}
\end{subfigure}
\hfill
\begin{subfigure}{0.32\textwidth}
\includegraphics[width=\textwidth]{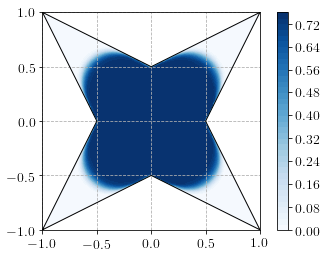}
\caption{$L^2$-norm}
\end{subfigure}
\hfill
\begin{subfigure}{0.32\textwidth}
\includegraphics[width=\textwidth]{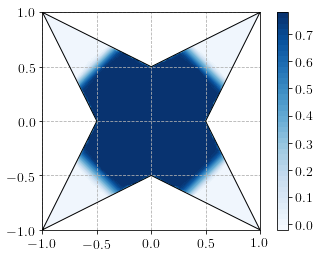}
\caption{$L^{\infty}$-norm}
\end{subfigure}
\end{center}
\caption{Generalized Cheeger sets of a star-shaped domain using different norms}
\label{Cheeger-star}
\end{figure}

These functions inherit from the composite class \texttt{ConvexFunction} which, by default, behaves like a \texttt{MeshConvexFunction}. To use them as \texttt{FacetConvexFunction}, they can be instantiated as \texttt{F = L2Norm.on\_facet(u)}. Using such predefined functions, many problems can be formulated in an extremely simple manner, without even worrying about the conic reformulation. For instance, we revisited the Cheeger problem on a star-shaped domain but with anisotropic norms \cite{kawohl2008p} such as $L^1$ and $L^{\infty}$ instead of $L^2$ in \eqref{Cheeger}\footnote{Note that in the general case of an $L^p$-norm for the gradient term, the corresponding jump term in \eqref{Cheeger-disc} is $\int_{\Gamma} |\jump{u}|\cdot\|\boldsymbol{n}\|_p\dS $ where $\boldsymbol{n}$ is the facet normal and similarly for the Dirichlet boundary term.}, the resulting sets are represented on Figure \ref{Cheeger-star}.

\section{A gallery of illustrative examples}\label{sec:gallery}
We now give a series of examples which illustrate the versatility of the \texttt{fenics\_optim} package for formulating and solving problems taken from the fields of solid and fluid mechanics, image processing and applied mathematics. The last two examples involve, in particular, time-dependent problems. Let us again point out that discretization choices or solver strategies using interior-point methods are not necessarily the most optimal ones for each of these problems and that many other approaches which have been proposed in the literature may be much more appropriate. We just aim at illustrating the potential of the package to formulate and solve various problems.

\subsection{Limit analysis of thin plates in bending}
The first problem consists in finding the ultimate load factor that a thin plate in bending can sustain given a predefined strength criterion and boundary condition. This limit analysis problem has been studied in \cite{demengel1983problemes, bleyer2016gamma}. In the present case, we consider a unit square plate made of a von Mises material of uniform bending strength $m$ and subjected to a uniformly distributed loading $f$. The thin plate limit analysis problem consists in solving the following problem:
\begin{equation}
\begin{array}{rl}
\displaystyle{\inf_{u\in \text{HB}_0(\Omega)}} & \displaystyle{\int_\Omega \pi(\nabla^2 u)\dx} \\
\text{s.t.} & \displaystyle{\int_{\Omega} fu \dx = 1}
\end{array} \label{BENDING-PLATE}
\end{equation}
where $\text{HB}_0$ is the space of bounded Hessian functions \cite{demengel1984fonctions} with zero trace on $\partial \Omega$ and $\pi(M) = \frac{2m}{\sqrt{3}}\sqrt{M_{11}^2+M_{22}^2+M_{12}^2+M_{11}M_{22}}$ for any $M\in \mathbb{S}_2^+$. One can notice that problem \eqref{BENDING-PLATE} shares some similar structure with the Cheeger problem \eqref{Cheeger} except that we are now dealing with the Hessian operator and a different norm through function $\pi$.

Contrary to elastic bending plate problems involving functions with $C^1$-continuity, we deal here with functions in HB which are continuous but may have discontinuities in their normal gradient $\partial_n u$, in particular we can consider again a Lagrange interpolation for $u$ with jumps of $\partial_n u $ across all internal facets $F\in \Gamma_h$ of unit normal $\boldsymbol{n}$. The $\pi$-function being some generalized total variation for $\nabla^2 u$, we have explicitly \cite{bleyer2013performance}:
\begin{equation}
\begin{array}{rl}
\displaystyle{\inf_{u\in \text{HB}_0(\Omega)}} & \displaystyle{\sum_{T\in \Tt_h} \int_T \pi(\nabla^2 u)\dx + \sum_{F\in \Gamma_h} \int_F \pi(\jump{\partial_n u} \boldsymbol{n} \otimes \boldsymbol{n})\dS} \\
\text{s.t.} & \displaystyle{\int_{\Omega} fu \dx = 1}
\end{array} 
\end{equation}
where it happens that in fact $\pi(\jump{\partial_n u} \boldsymbol{n} \otimes \boldsymbol{n}) = |\jump{\partial_n u}|\pi(\boldsymbol{n}\otimes \boldsymbol{n}) = |\jump{\partial_n u}|\frac{2m}{\sqrt{3}}$. Following \ref{appendix:plate}, we have the following formulation of the bending plate problem for a $\PP^2$ interpolation:
\begin{pythoncode}
prob = MosekProblem("Bending plate limit analysis")

V = FunctionSpace(mesh, "CG", 2)
bc = DirichletBC(V, Constant(0.), boundary)
u = prob.add_var(V, bc = bc)

R = FunctionSpace(mesh, "R", 0)
def Pext(lamb):
    return [lamb*dot(load,u)*dx]
prob.add_eq_constraint(R, A=Pext, b=1)

J = as_matrix([[2., 1., 0.],
               [0, sqrt(3.), 0.],
               [0, 0, 1]]) 
def Chi(v):
    chi = sym(grad(grad(v)))
    return as_vector([chi[0,0], chi[1,1], 2*chi[0, 1]])
pi_c = L2Norm(u, "vertex", degree=1)
pi_c.set_term(m/sqrt(3)*dot(J, Chi(u)))
prob.add_convex_term(pi_c)

pi_h = L1Norm.on_facet(u)
pi_h.set_term([jump(grad(u), n)], k=2/sqrt(3)*m)
prob.add_convex_term(pi_h)

prob.optimize()
\end{pythoncode}

\begin{figure}
\begin{center}
\begin{subfigure}{0.49\textwidth}
\includegraphics[width=\textwidth]{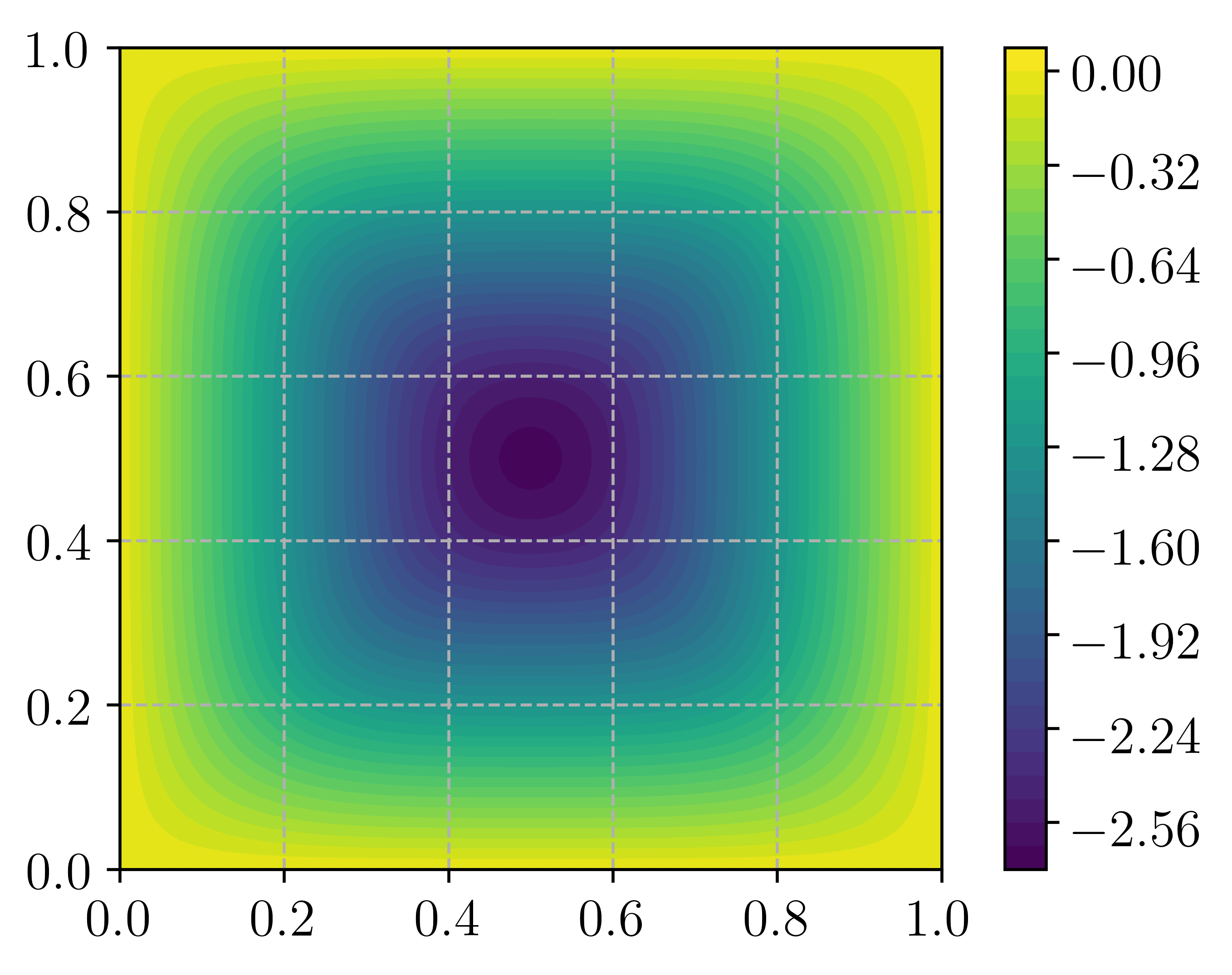}
\caption{Optimal collapse mechanism $u$}
\end{subfigure}
\hfill
\begin{subfigure}{0.49\textwidth}
\includegraphics[width=0.95\textwidth]{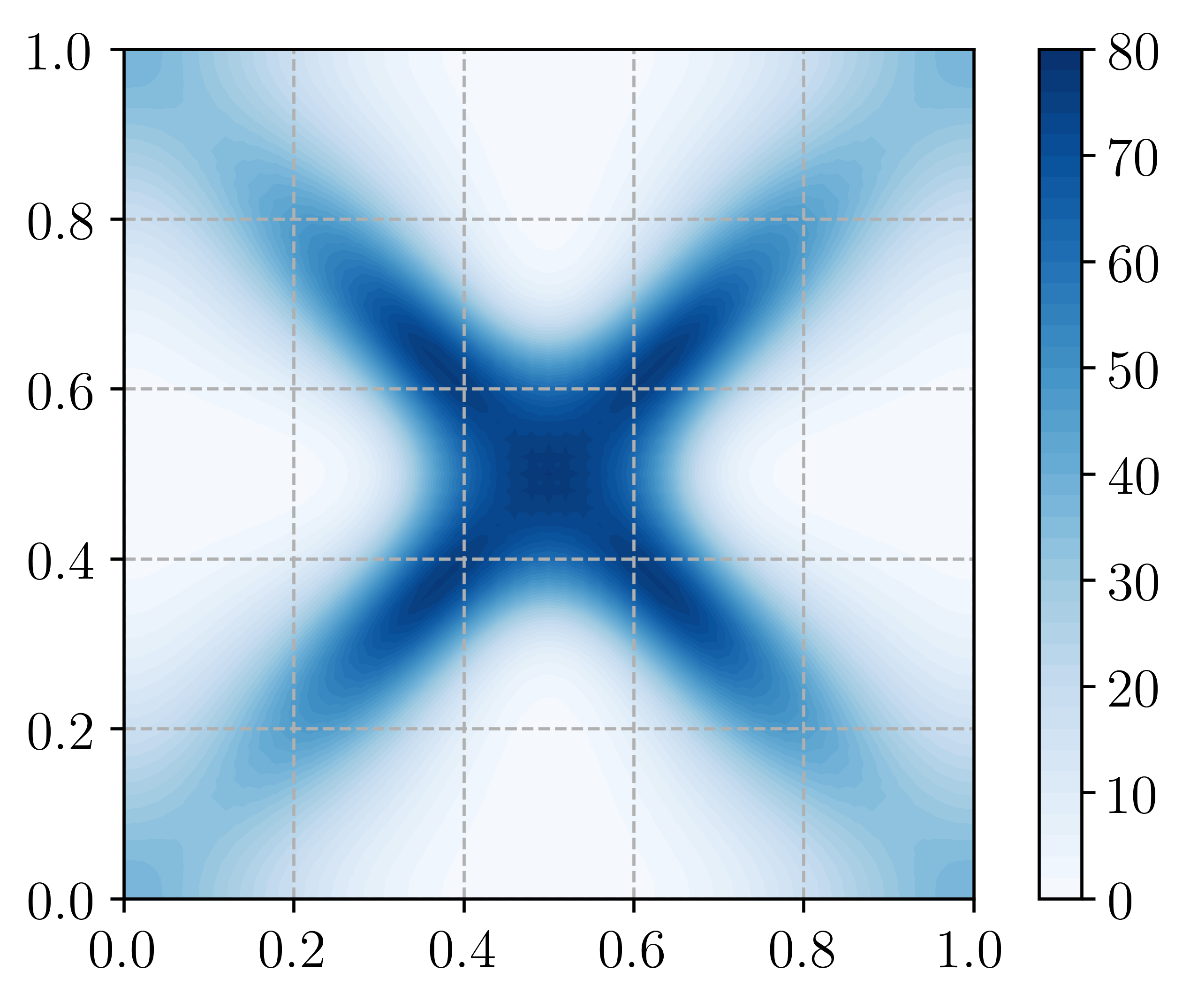}
\caption{Curvature dissipation density $\pi(\nabla^2 u)$}
\end{subfigure}
\end{center}
\caption{Results for the simply supported von Mises square plate}
\label{plate-results}
\end{figure}

The reference solution for this problem is known to be $25.02 m/f$ \cite{capsoni1999limit}, whereas we find $25.05 m/f$ for a $50\times 50$ structured mesh. The corresponding solutions for $u$ and $\pi(\nabla^2 u)$ are represented in Figure \ref{plate-results}.

\subsection{Viscoplastic yield stress fluids}
Viscoplastic (or yield stress) fluids \cite{balmforth2014yielding,coussot2016bingham} are a particular class of non-Newtonian fluids which, in their most simple form, namely the Bingham model, behave like a purely rigid solid when the shear stress is below a critical yield stress $\tau_0$ and flow like a Newtonian fluid when the shear stress is above $\tau_0$. They appear in many applications ranging from civil engineering, petroleum, cosmetics or food industries. The solution of a steady state viscoplastic fluid flow under Dirichlet boundary conditions and a given external force field $\boldsymbol{f}$ can be obtained as the unique solution to the following convex variational principle \cite{glowinski1984}:
\begin{equation}
\begin{array}{rl}
\displaystyle{\inf_{\boldsymbol{u}\in \boldsymbol{H}^1(\Omega;\RR^d)}} & \displaystyle{\int_\Omega (\mu\|\nabla \boldsymbol{u}\|_2^2 + \sqrt{2}\tau_0\|\nabla \boldsymbol{u}\|_2)\dx - \int_{\Omega}\boldsymbol{f}\cdot\boldsymbol{u}\dx} \\
\text{s.t.} & \div \boldsymbol{u} = 0 \text{ in }\Omega \\
& \boldsymbol{u} = \boldsymbol{g} \text{ on } \partial \Omega
\end{array} \label{viscoplastic-fluids}
\end{equation}
where $\mu$ is the fluid viscosity. Typical solutions of problem \eqref{viscoplastic-fluids} involve rigid zones in which $\nabla \boldsymbol{u} =0$ and flowing regions where $\|\nabla\boldsymbol{u}\|\neq 0$, the locations of which are \textit{a priori} unknown. Note that when $\tau_0=0$, we recover the classical viscous energy of Stokes flows and optimality conditions of problem \eqref{viscoplastic-fluids} reduce to a linear problem. The FE discretization is quite classical, we adopt Taylor-Hood $\PP^2/\PP^1$ discretization for the velocity $\boldsymbol{u}$ and the pressure $p$ which is the Lagrange multiplier of constraint $\div \boldsymbol{u}=0$.

The considered problem is the classical lid-driven unit-square cavity, with $\boldsymbol{f}=0$, $\boldsymbol{u}=0$ everywhere on $\partial \Omega$, except on the top boundary $y=1$ where $\boldsymbol{u}=(U, 0)$ with $U$ the imposed constant velocity. Different solutions to problem \eqref{viscoplastic-fluids} are then obtained depending on the value of the non-dimensional Bingham number $\texttt{Bi}=\dfrac{\mu U}{\tau_0 L}$ with the characteristic length $L=1$ for the present case. When $\texttt{Bi}=0$, the solution is that of a Newtonian fluid and when $\texttt{Bi}\to\infty$ it corresponds to that of a purely plastic material.

Implementation in \texttt{fenics\_optim} is straightforward once the symmetric tensor $\nabla\boldsymbol{u}$ has been represented as a vector of $\RR^3$ through the \texttt{strain} function \cite{bleyer2018advances}.
\begin{pythoncode}
prob = MosekProblem("Viscoplastic fluid")

V = VectorFunctionSpace(mesh, "CG", 2)
bc = [DirichletBC(V, Constant((1.,0.)), top),
      DirichletBC(V, Constant((0.,0.)), sides)]
u = prob.add_var(V, bc=bc)

Vp = FunctionSpace(mesh, "CG", 1)
def mass_conserv(p):
    return [p*div(u)*dx]
prob.add_eq_constraint(Vp, mass_conserv)

def strain(v):
    E = sym(grad(v))
    return as_vector([E[0, 0], E[1, 1], sqrt(2)*E[0, 1]])
visc = QuadraticTerm(u, degree=2)
visc.set_term(strain(u))
plast = L2Norm(u, degree=2)
plast.set_term(strain(u))

prob.add_convex_term(2*mu*visc)
prob.add_convex_term(sqrt(2)*tau0*plast)

prob.optimize()
\end{pythoncode}

\begin{figure}
\begin{center}
\begin{subfigure}{0.49\textwidth}
\includegraphics[width=\textwidth]{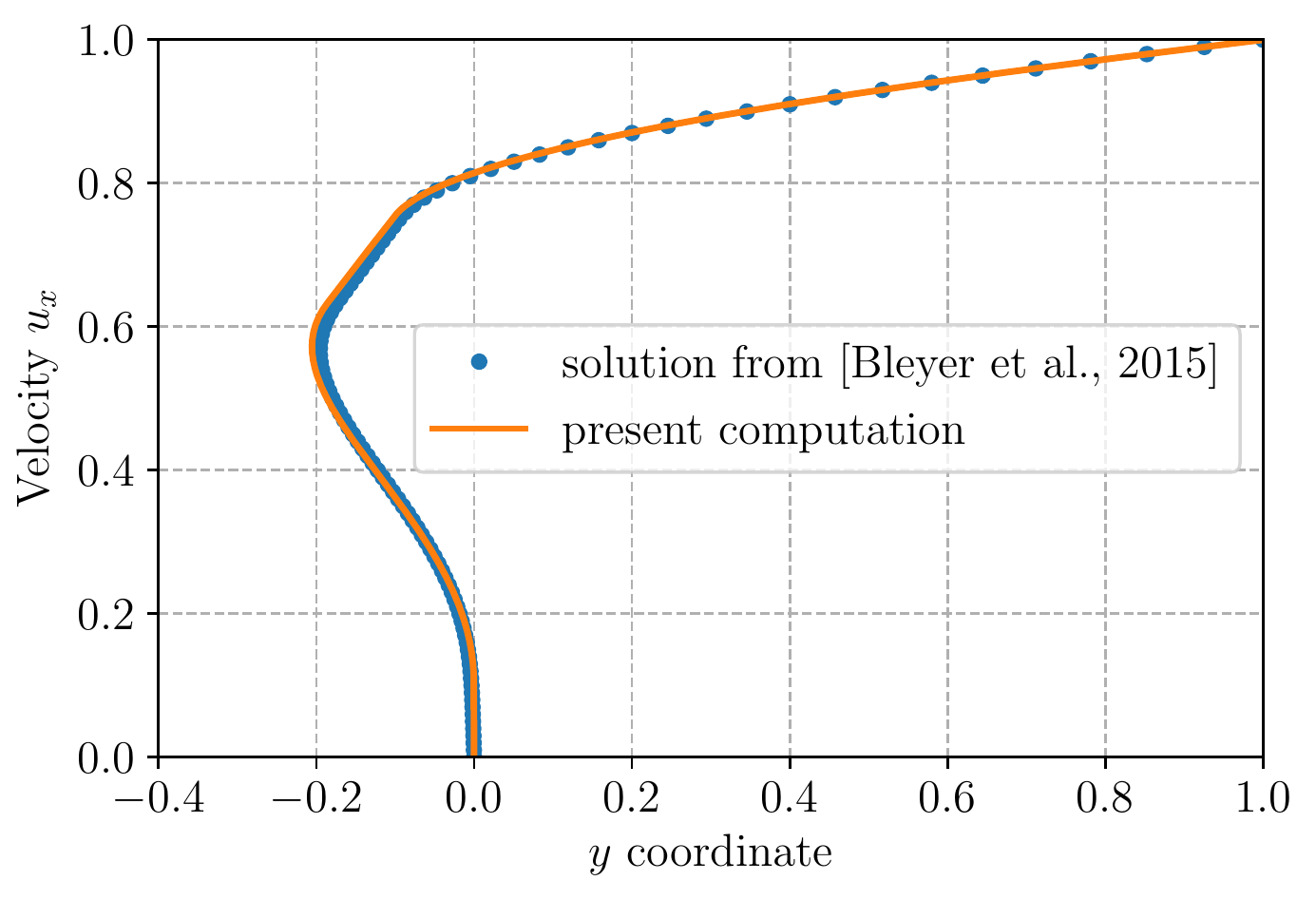}
\caption{$\texttt{Bi}=2$}
\end{subfigure}
\hfill
\begin{subfigure}{0.49\textwidth}
\includegraphics[width=\textwidth]{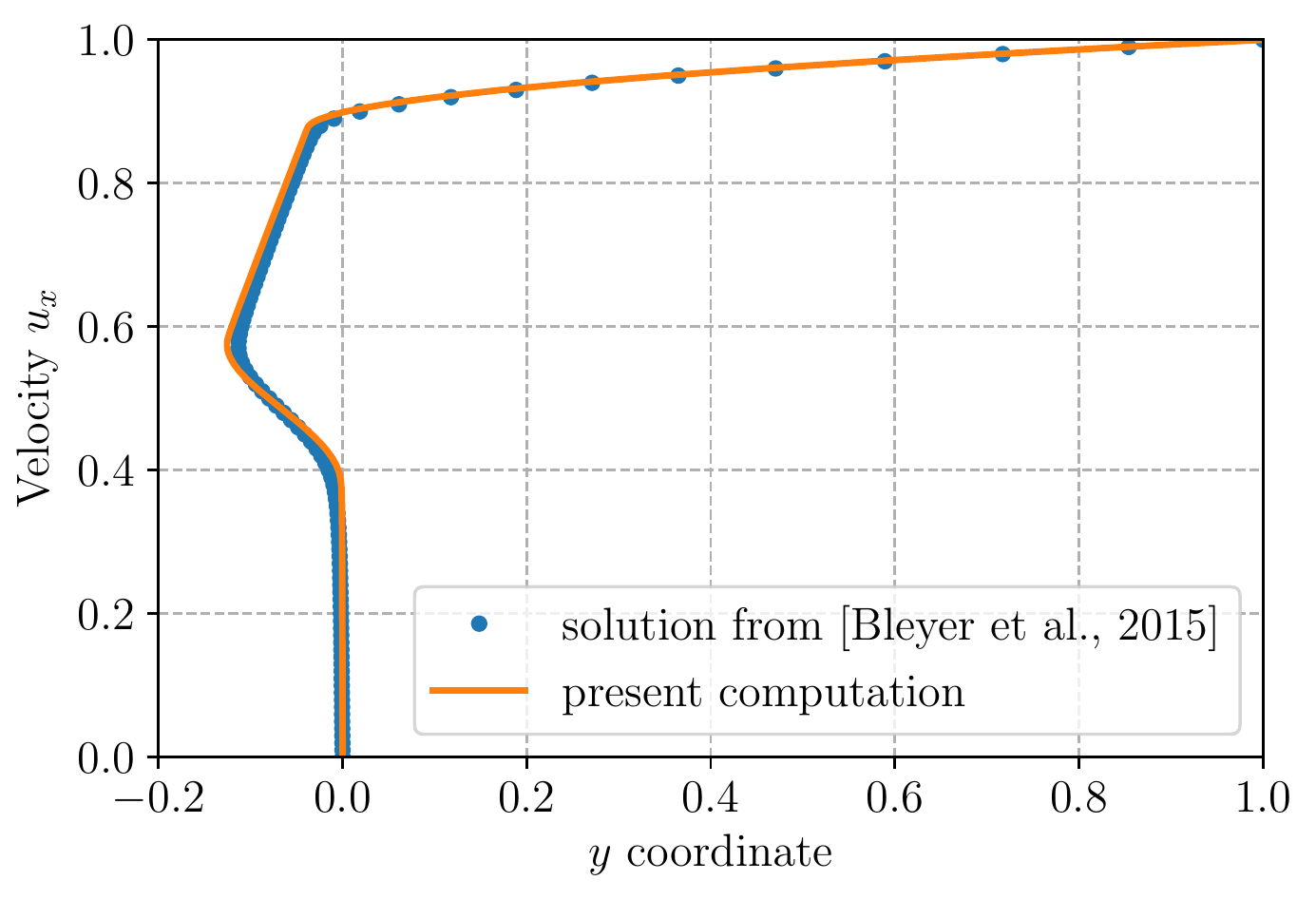}
\caption{$\texttt{Bi}=20$}
\end{subfigure}
\end{center}
\caption{Horizontal velocity profile $u_x(y)$ on the middle plane $x=0.5$, comparison with results from \cite{bleyer2015efficient}}
\label{viscoplastic-profile}
\end{figure}

\begin{figure}
\begin{center}
\begin{subfigure}{0.49\textwidth}
\includegraphics[width=\textwidth]{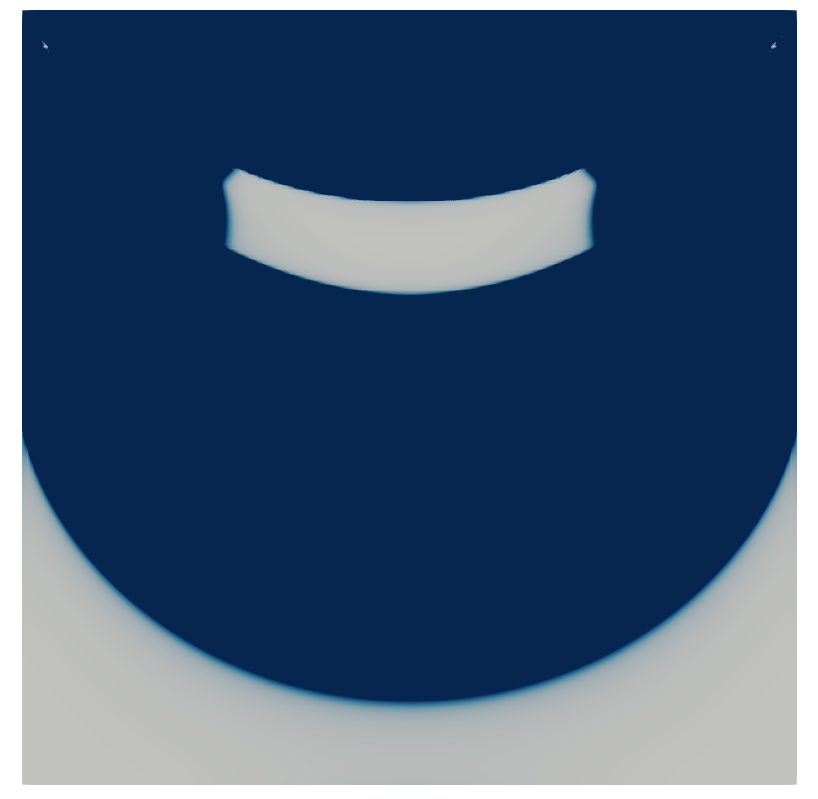}
\caption{$\texttt{Bi}=2$}
\end{subfigure}
\hfill
\begin{subfigure}{0.49\textwidth}
\includegraphics[width=\textwidth]{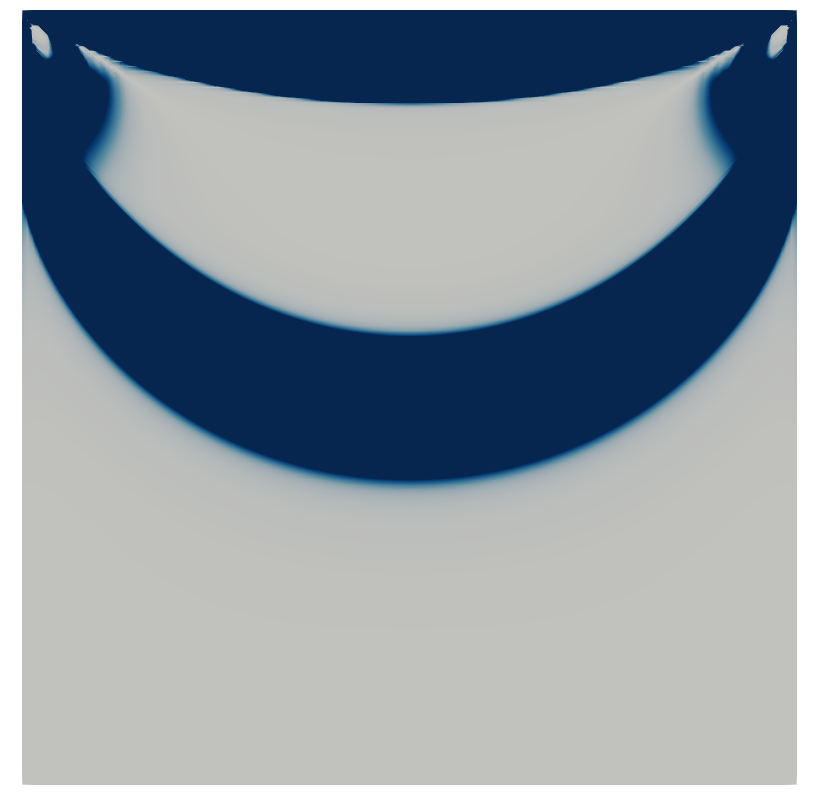}
\caption{$\texttt{Bi}=20$}
\end{subfigure}
\end{center}
\caption{Transition between solid (white) and liquid (regions). The bottom solid region is arrested and the central region rotates like a rigid body.}
\label{yield-surfaces}
\end{figure}
The obtained optimal velocity field is compared on Figure \ref{viscoplastic-profile} with that from a previous independent implementation described in \cite{bleyer2015efficient}. Finally, if $\boldsymbol{d}=\nabla \boldsymbol{u}\neq 0$, then the stress inside the fluid is given by $\boldsymbol{\tau}=2\mu\boldsymbol{d}+\sqrt{2}\tau_0\dfrac{\boldsymbol{d}}{\|\boldsymbol{d}\|_2}$ and $\|\boldsymbol{\tau}\|_2 > \sqrt{2}\tau_0$. In Figure \ref{yield-surfaces}, $\|\boldsymbol{\tau}\|_2$ has been plotted with a colormap ranging from $\sqrt{2}\tau_0$ to $1.01\sqrt{2}\tau_0$, thus exhibiting the transition from solid regions (white) to liquid regions (blue).

\subsection{Total Variation inpainting}
\begin{figure}
\begin{center}
\begin{subfigure}{\textwidth}
\includegraphics[width=\textwidth]{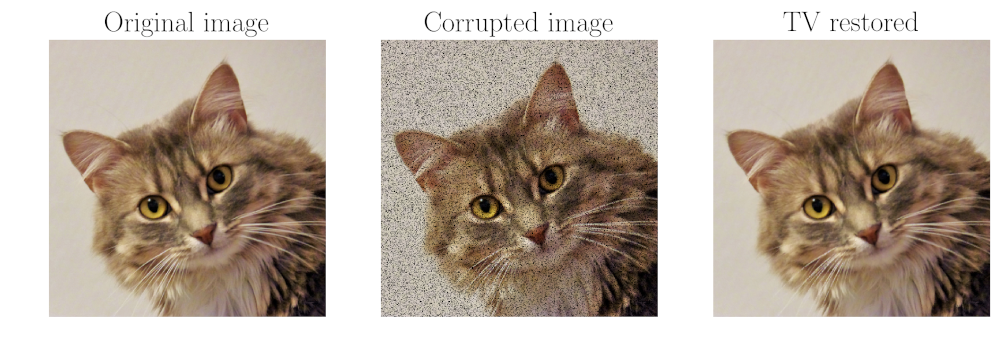}
\caption{$\eta=25\%$ corruption level}
\end{subfigure}
\hfill
\begin{subfigure}{\textwidth}
\includegraphics[width=\textwidth]{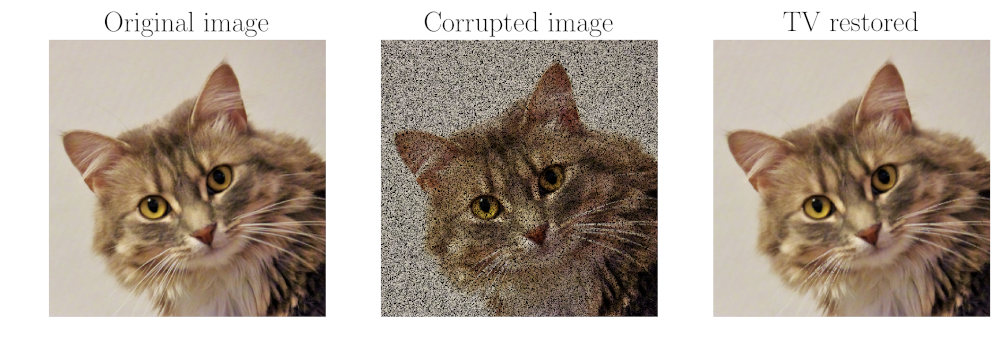}
\caption{$\eta=50\%$ corruption level}
\end{subfigure}
\end{center}
\caption{Inpainting problem of a corrupted image using TV restoration}
\label{inpainting-results}
\end{figure}

In this example, we consider an image processing problem called \textit{inpainting}, consisting in recovering an image which has been deteriorated. In the present case, we consider a color RGB image in which a fraction $\eta$ of randomly chosen pixels have been lost (black). The inpainting problem consists in recovering the three color channels $\mathbf{U} = (\mathbf{u}_{j})$ for $j\in \{R,G,B\}$ such that it matches the original color for pixels which have not been corrupted and minimizing a given energy for the remaining pixels. An efficient choice of energy for the inpainting problem is the $L^2$ total variation norm $TV(u) = \int_\Omega \|\nabla u\|_2\dx$ for a given color channel $u$. For an image, the discrete gradient can be computed by finite differences. Here, as we work with a FE library, the image will be represented using a Crouzeix-Raviart ($CR$) interpolation \cite{chambolle2018crouzeix} on a structured finite element mesh. The inpainting problem therefore reads as:

\begin{equation}
\begin{array}{rl} \displaystyle{\min_{\mathbf{U}\in(\RR^N)^3}} & \displaystyle{\sum_{j\in\{R,G,B\}} \int_\Omega \|\nabla u_j\|_2 \dx}  \\
\text{s.t.}  & U_{i,j} = U_{i,j}^\text{orig} \quad \forall i\notin I_c,\:\: \forall j\in\{R,G,B\}
\end{array} \label{inpainting}
\end{equation}
where $I_c$ denotes the set of corrupted pixels. Again, problem \eqref{inpainting} can be defined very easily as follows:
\begin{pythoncode}
prob = MosekProblem("TV inpainting")
u = prob.add_var(V, ux=ux, lx=lx)

for i in range(3):
    tv_norm = L2Norm(u)
    tv_norm.set_term(grad(u[i]))
    prob.add_convex_term(tv_norm)

prob.optimize()
\end{pythoncode}
where \texttt{V} is the space $(CR)^3$ and \texttt{ux} (resp. \texttt{lx}) denote functions of \texttt{V} equal to the original image on cells corresponding to uncorrupted pixels and which take $+\infty$ (resp. $-\infty$) values on $I_c$, so that $\texttt{lx} \leq \texttt{u} \leq \texttt{ux}$ amounts to enforcing fidelity with the uncorrupted values. Finally, an $L^2$-norm term on the gradient of each channel is added to the problem. Results for a 512$\times$512 image discretized using a triangular mesh of identical resolution (each pixel is split into two triangles) are represented in Figure \ref{inpainting-results} for two corruption levels. It must be noted that optimization took roughly one minute for both cases.

\subsection{Cartoon+Texture Variational Image Decomposition}
The next image processing example we consider is that of decomposing an image $y=u+v$ into a cartoon-like component $u$ and a texture component $v$ (here we assume that the image is not noisy). The cartoon layer $u$ captures flat regions separated by sharp edges, whereas the texture component $v$ contains the high frequency oscillations. There are many existing models to perform such a decomposition, in the following, we implement the model proposed by Y. Meyer \cite{meyer2001oscillating, weiss2009efficient}:
\begin{equation}
\begin{array}{rl}
\displaystyle{\inf_{u,v}} & \displaystyle{\int_\Omega \|\nabla u\|_2 \dx + \alpha \|v\|_G} \\
\text{s.t.} & y=u+v \\
\\
\text{where} &  \displaystyle{\|v\|_G = \inf_{\boldsymbol{g}\in L^{\infty}(\Omega;\RR^2)} \{\|\sqrt{g_1^2+g_2^2}\|_{\infty} \text{ s.t. } v=\div \boldsymbol{g}\}}
\end{array} 
\end{equation}
This model favors flat regions in $u$ due to the use of the TV norm and oscillatory regions in $v$ since $\|v\|_G$ increases for characteristic functions. Following \cite{weiss2009efficient}, we reformulate the model as:
\begin{equation}
\begin{array}{rl}
\displaystyle{\inf_{u,\boldsymbol{g}}} & \displaystyle{\int_\Omega \|\nabla u\|_2 \dx} \\
\text{s.t.} & y=u+\div(\boldsymbol{g})
\\ & \sqrt{g_1^2+g_2^2}\leq \alpha
\end{array} 
\end{equation}
The original image (512$\times$512) is here represented on a triangular finite-element mesh of similar mesh size and we adopt a Crouzeix-Raviart interpolation for $u$ and a Raviart-Thomas interpolation for $\boldsymbol{g}$. The constraint $ y=u+\div(\boldsymbol{g})$ is enforced weakly on the CR space. The implementation reads as:
\begin{pythoncode}
prob = MosekProblem("Cartoon/texture decomposition")
Vu = FunctionSpace(mesh, "CR", 1)
Vg = FunctionSpace(mesh, "RT", 1)
u, g = prob.add_var([Vu, Vg])

def constraint(l):
    return [dot(l, u)*dx, dot(l, div(g))*dx]
def rhs(l):
    return dot(l, y)*dx
prob.add_eq_constraint(Vu, A=constraint, b=rhs)

tv_norm = L2Norm(u)
tv_norm.set_term(grad(u))
prob.add_convex_term(tv_norm)
    
g_norm = L2Ball(g)
g_norm.set_term(g, k=alpha)
prob.add_convex_term(g_norm)

prob.optimize()
\end{pythoncode}

\begin{figure}
\begin{center}
\includegraphics[width=\textwidth]{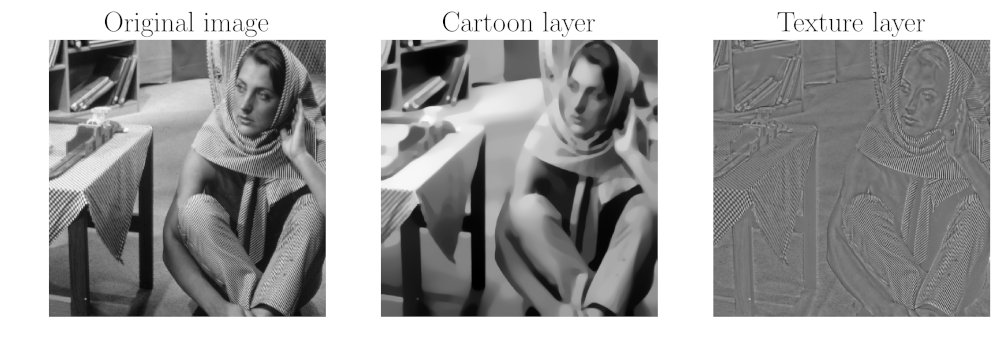}
\end{center}
\caption{Cartoon-texture decomposition with $\alpha=2\text{e-}4$}
\label{barbara}
\end{figure}

Results for the Barbara image decomposition are shown in Figure \ref{barbara}.

\subsection{Time-dependent sandpile growth}
\begin{figure}
\begin{center}
\begin{subfigure}{0.49\textwidth}
\includegraphics[width=\textwidth]{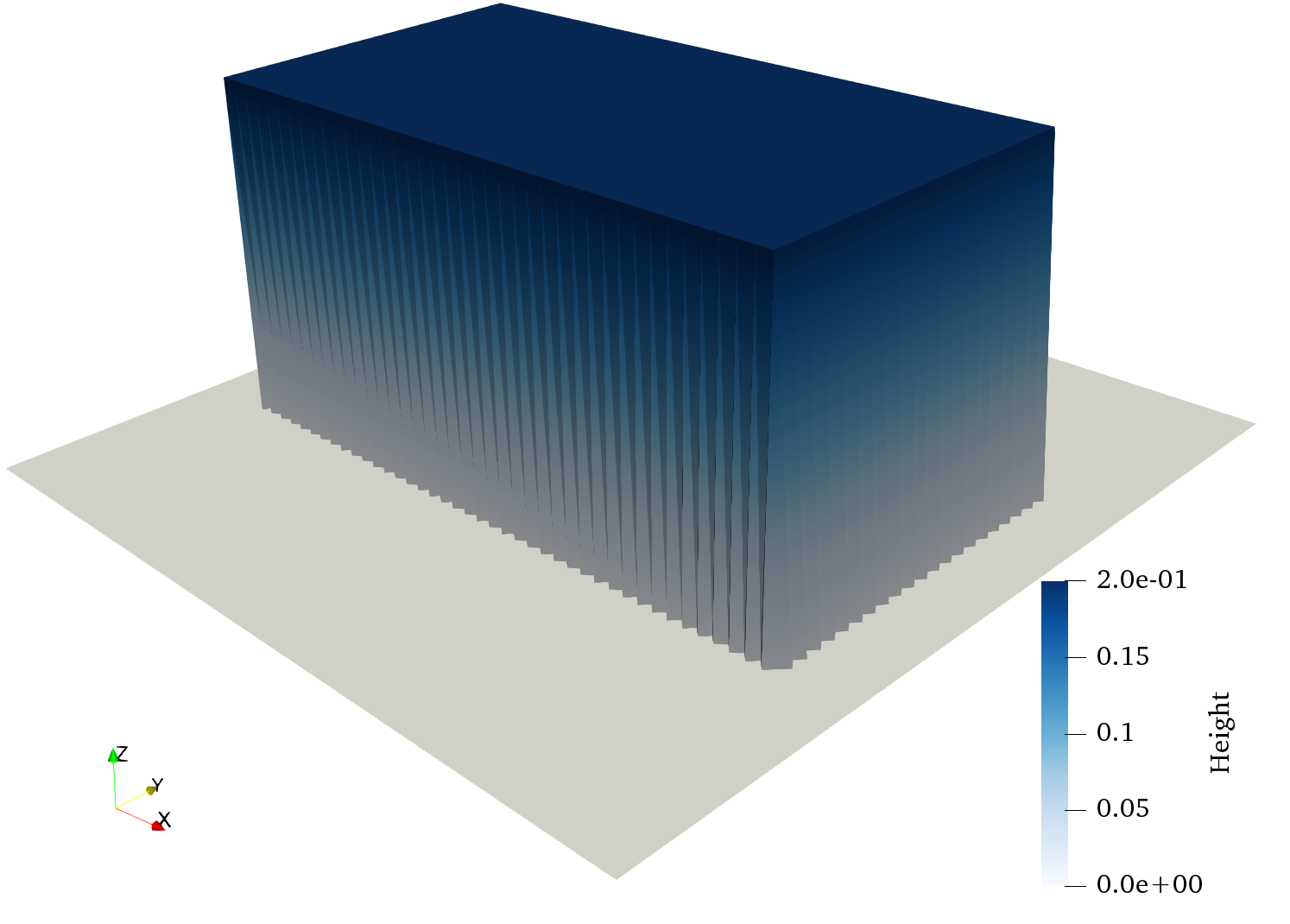}
\caption{$t=0$}
\end{subfigure}
\hfill
\begin{subfigure}{0.49\textwidth}
\includegraphics[width=\textwidth]{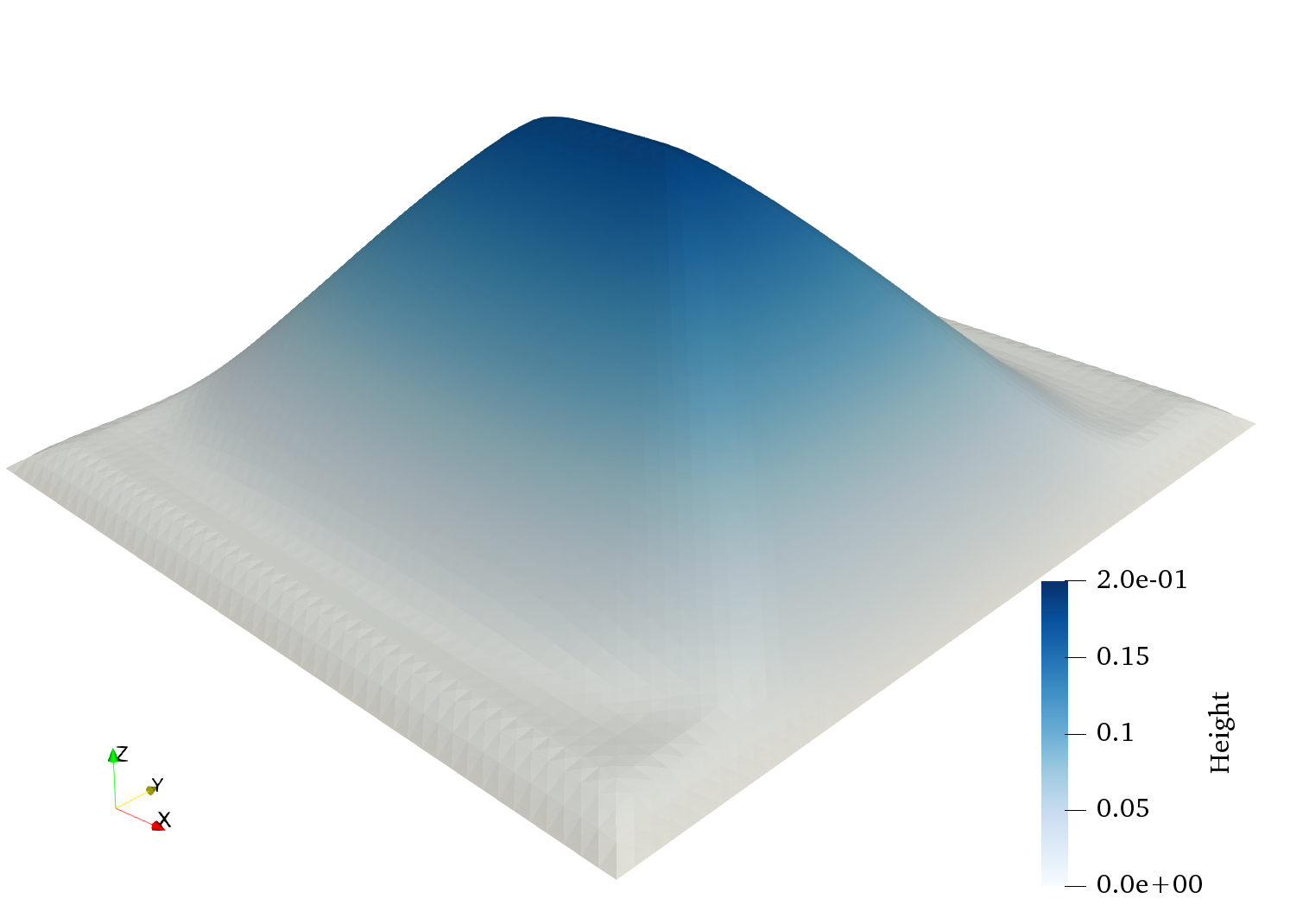}
\caption{$t=0.05$}
\end{subfigure}\\
\begin{subfigure}{0.49\textwidth}
\includegraphics[width=\textwidth]{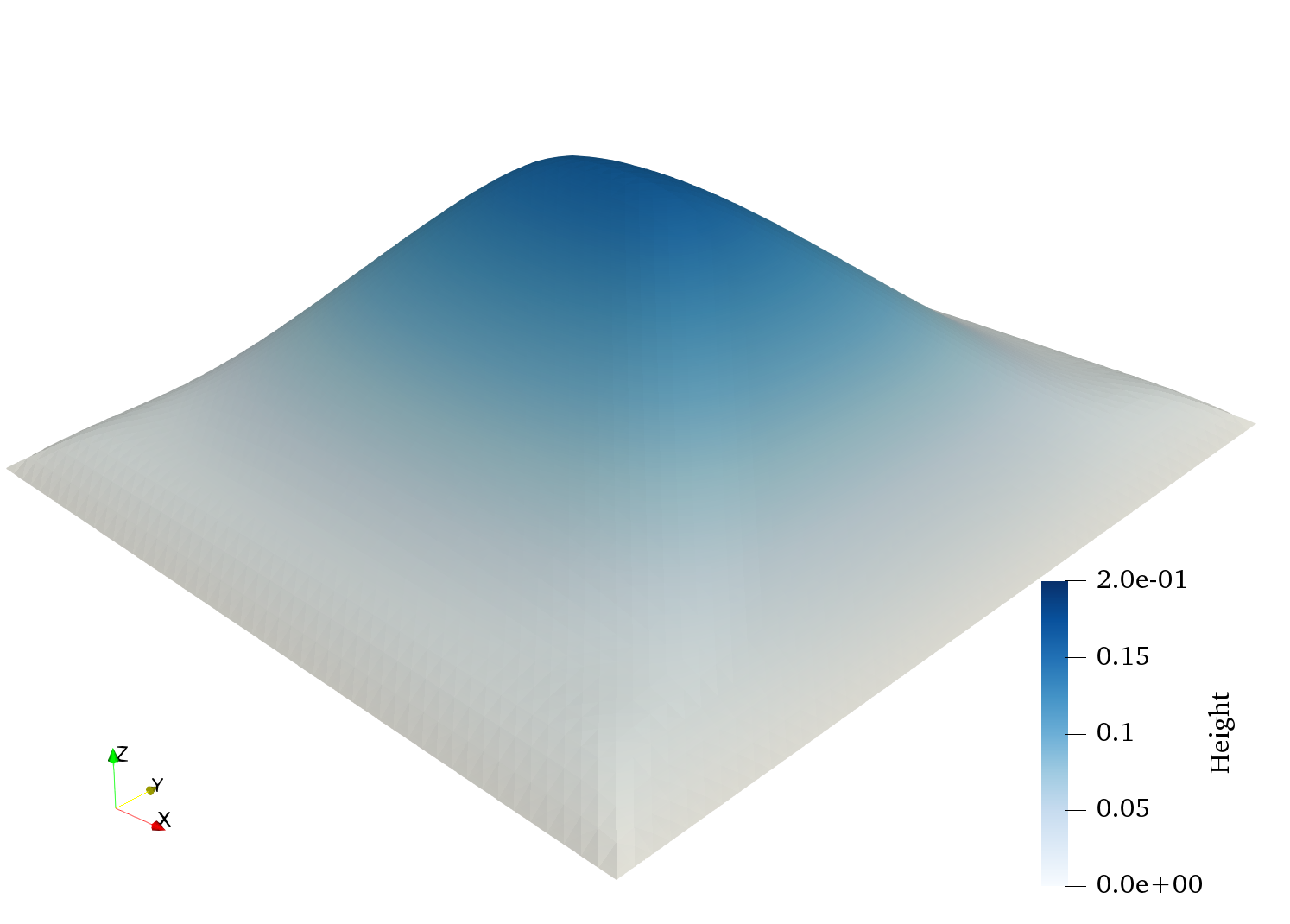}
\caption{$t=0.1$}
\end{subfigure}
\hfill
\begin{subfigure}{0.49\textwidth}
\includegraphics[width=\textwidth]{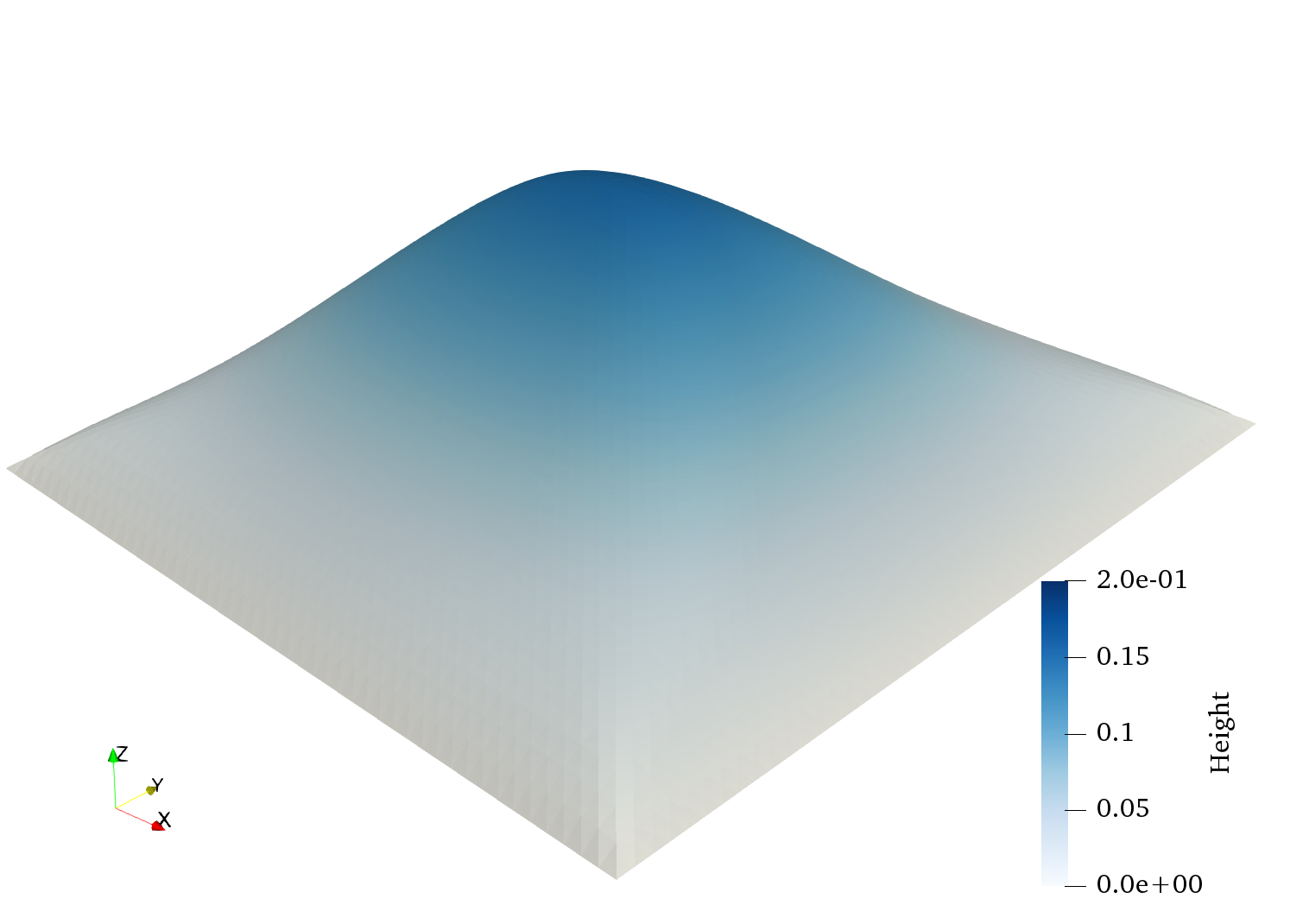}
\caption{$t=0.25$}
\end{subfigure}
\end{center}
\caption{Sandpile growth evolution starting from an initial unstable configuration (height amplification by factor 2)}
\label{sandpile-collapse}
\end{figure}

In this example, we consider the time-dependent evolution model of a sandpile characterized by its height $h$. Since sand can fall off the table domain $\Omega$, Dirichlet boundary conditions are prescribed. Layers of sand having a slope larger than the critical angle at rest $\tan\alpha$ will fall down the slope and can be modelled by Prighozin evolutionary PDE \cite{prigozhin1996variational}:
\begin{align}
\partial_t h - \div(m\nabla h) &= f \quad \text{in } \Omega\times[0;T] \\
m\geq 0, \: \|\nabla h\|\leq \tan \alpha & \notag \\ 
m( \|\nabla h\|-\tan \alpha) &=0 \notag
\end{align}
with $h(\boldsymbol{x}, t)=0$ for $\boldsymbol{x}\in\partial \Omega$ and $h(\boldsymbol{x},0) = h_0(\boldsymbol{x})$ as the initial sandpile height. In the above, $-m\nabla h$ denotes the horizontal material flux of collapsing sand layers and $f$ is a potential source term. This model has been also linked with the Monge-Kantorovitch problem of optimal mass transportation. Performing a backward implicit Euler discretization of the time derivative at each time step $t_n$ and knowing the previous height configuration $h_{n-1}(\boldsymbol{x})$ at time $t=t_{n-1}$, finding $h_n(\boldsymbol{x})$ amounts to solving the following variational problem \cite{dumont2009dual}:
\begin{equation}
\begin{array}{rl}
\displaystyle{\inf_{h}} & \displaystyle{\dfrac{1}{2}\int_\Omega (h-g_n)^2\dx} \\
\text{s.t.} & \|\nabla h\|\leq \tan \alpha 
\end{array} \label{sandpile}
\end{equation}
where $g_n=\Delta t f+h_{n-1}$ and $\Delta t=T/N$ is the time interval of each $N$ time increments discretization of interval $[0;T]$. Adopting a standard Lagrange $\PP^1$ interpolation for $h$, problem \eqref{sandpile} is solved $N$ times with values of $g_n$ updated from the previous solution. Figure \ref{sandpile-collapse} illustrate the results obtained with $\alpha=30^{\circ}$, no source term $f=0$ and an initial unstable configuration for $h_0$ since $\|\nabla h_0\| > \tan \alpha$.

\subsection{Optimal transport with space-time finite elements}
Finally, we consider the Brenier-Benamou dynamic formulation \cite{benamou2000computational} of quadratic cost optimal transport between two distributions $\rho_0$ and $\rho_1$ which reads as:
\begin{equation}
\begin{array}{rl}
\displaystyle{\inf_{\rho,\boldsymbol{v}}} & \displaystyle{\dfrac{1}{2}\int_0^1 \int_\Omega \rho(\boldsymbol{x},t)\|\boldsymbol{v}(\boldsymbol{x},t)\|_2^2\dx\dt} \\
\text{s.t.} & \partial_t \rho + \div_x(\rho\boldsymbol{v})=0 \\
& \rho(\boldsymbol{x},t=0)=\rho_0(\boldsymbol{x}) \\
& \rho(\boldsymbol{x},t=1)=\rho_1(\boldsymbol{x}) \\
& \boldsymbol{v}\cdot\boldsymbol{n} = 0 \text{ on }\partial \Omega
\end{array} \label{optimal-transport}
\end{equation}

The change of variable $(\rho,\boldsymbol{m}) := (\rho,\rho\boldsymbol{v})$ proposed in \cite{benamou2000computational} enables to obtain the following convex optimization problem:
\begin{equation}
\begin{array}{rl}
\displaystyle{\inf_{\rho,\boldsymbol{m}}} & \displaystyle{\int_0^1 \int_\Omega c(\rho,\boldsymbol{m})\dx\dt} \\
\text{s.t.} & \partial_t \rho + \div_x\boldsymbol{m}=0 \\
& \rho(\boldsymbol{x},t=0)=\rho_0(\boldsymbol{x}) \\
& \rho(\boldsymbol{x},t=1)=\rho_1(\boldsymbol{x}) \\
& \boldsymbol{m}\cdot\boldsymbol{n} = 0 \text{ on }\partial \Omega
\end{array} \label{optimal-transport-cvx}
\end{equation}
where the cost function is $ c(\rho,\boldsymbol{m}) = \begin{cases} \dfrac{\|\boldsymbol{m}\|_2^2}{2\rho} & \text{if }\rho > 0 \\
0 & \text{if } (\rho,\boldsymbol{m})=(0,0) \\
+\infty & \text{otherwise} \end{cases}$. This function is convex and, observing that $c(\rho,\boldsymbol{m}) \leq t$ is equivalent to $2\rho t \geq \|\boldsymbol{m}\|^2_2$, is conic representable as follows:
\begin{equation}
\begin{array}{rl}
c(\rho,\boldsymbol{m}) = \displaystyle{\min_\mathbf{y}} & y_0 \\
\text{s.t.} & \begin{bmatrix}
\rho \\ m_1 \\ m_2
\end{bmatrix} - \begin{bmatrix}
y_1 \\ y_2 \\ y_3
\end{bmatrix} = 0 \\
& \mathbf{y}\in\Qq_4^r
\end{array} \label{cost-reformulation}
\end{equation}

The numerical approximation is performed by relying on a space-time finite element discretization of $Q=[0;1]^2\times[0;T]$ with $T=1$. We adopt $\PP^2$ Lagrange finite elements for the 3d-vector $(\rho,\boldsymbol{m})$. Initial and boundary conditions are imposed on the different boundaries of the space-time cube. The mass conservation equation is replaced by a relaxed inequality between $\pm\epsilon$ to allow for small deviations because of errors induced by the space-time discretization. It is written as:
\begin{pythoncode}
prob = MosekProblem("Optimal transport")
u = prob.add_var(V, bc=bc)

conserv = InequalityConstraint(u, degree=2)
rho, mx, my = u[0], u[1], u[2]
eps = 1e-6
conserv.set_term(rho.dx(2)+mx.dx(0)+my.dx(1), bl=-eps, bu=eps)
prob.add_convex_term(conserv)
\end{pythoncode}
where \texttt{dx(0)} and \texttt{dx(1)} stand for derivation along both spatial directions and \texttt{dx(2)} stands for derivation along the third time direction. Finally, the cost function term is added following reformulation \eqref{cost-reformulation}:
\begin{pythoncode}
class CostFunction(ConvexFunction):
    def conic_repr(self, X):
        Y = self.add_var(dim=4, cone=RQuad(4))
        Ybar = as_vector([Y[i] for i in range(1, 4)])
        self.add_eq_constraint([X, -Ybar])
        self.set_linear_term([None, Y[0]])
        
c = CostFunction(u, degree=2)
c.set_term(u)
prob.add_convex_term(c)
\end{pythoncode}

\begin{figure}
\begin{center}
\begin{subfigure}{0.32\textwidth}
\includegraphics[width=\textwidth]{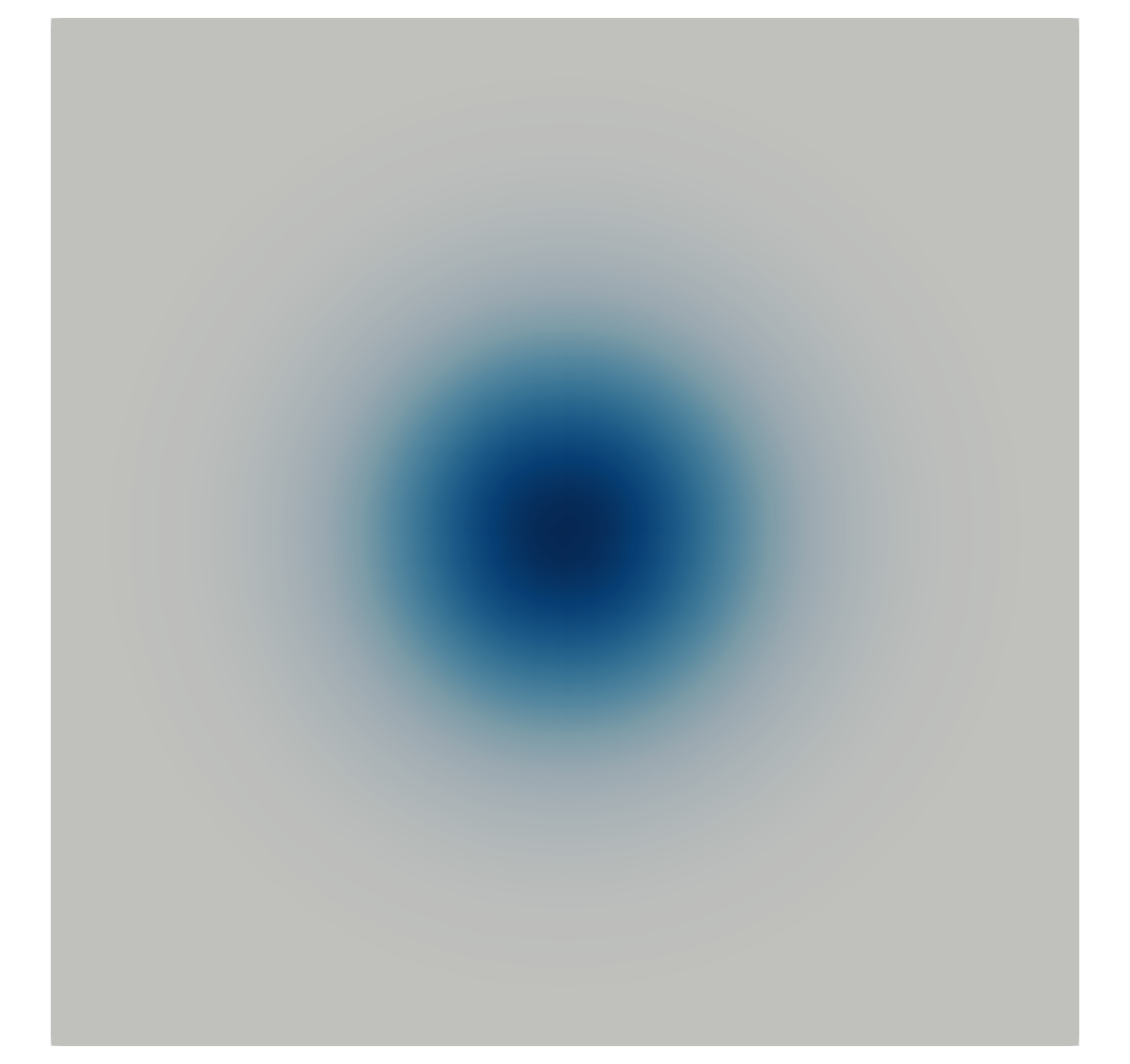}
\caption{$t=0$: $\rho_0$ distribution}
\end{subfigure}
\hfill
\begin{subfigure}{0.32\textwidth}
\includegraphics[width=\textwidth]{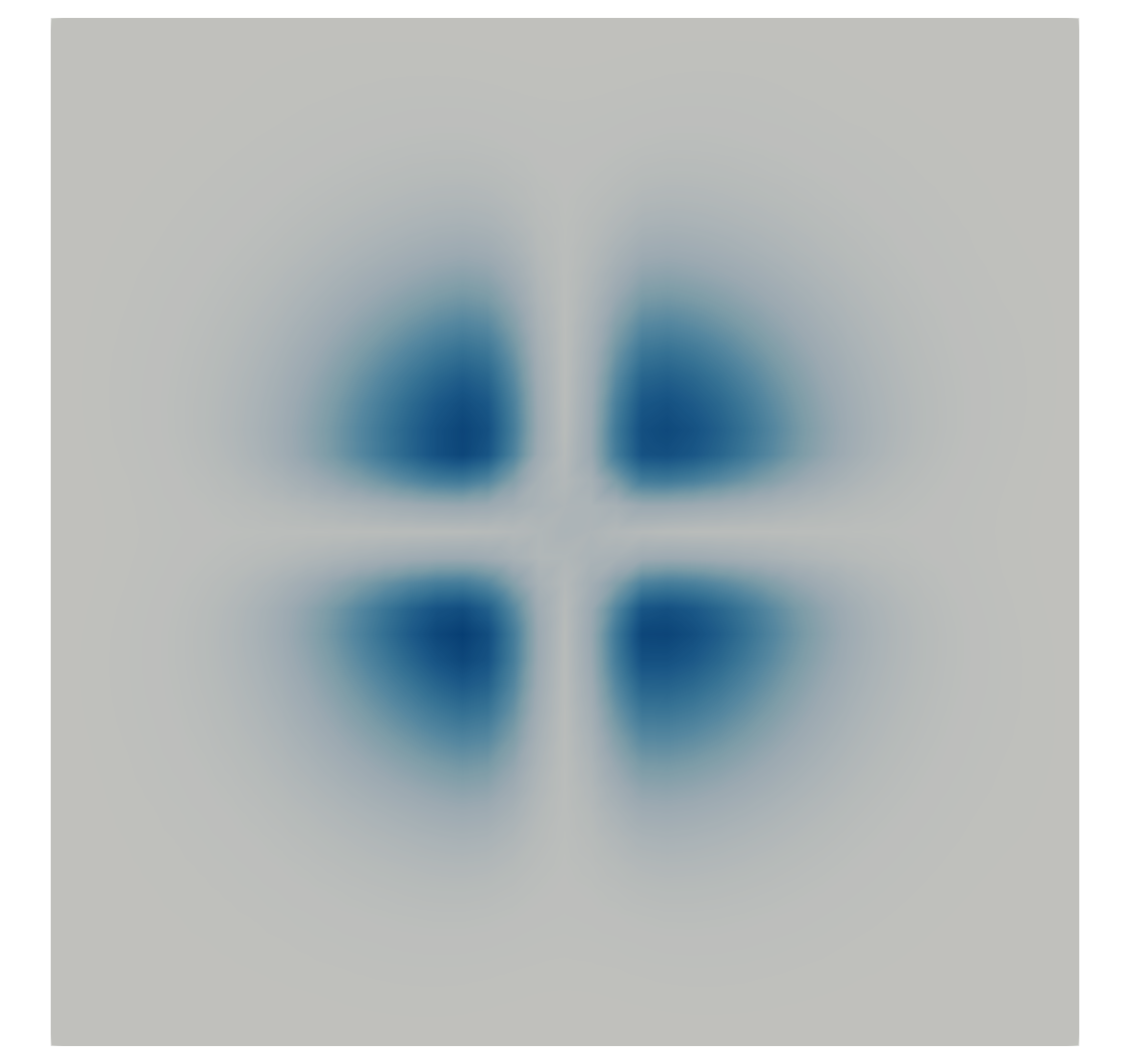}
\caption{$t=0.2$}
\end{subfigure}
\hfill
\begin{subfigure}{0.32\textwidth}
\includegraphics[width=\textwidth]{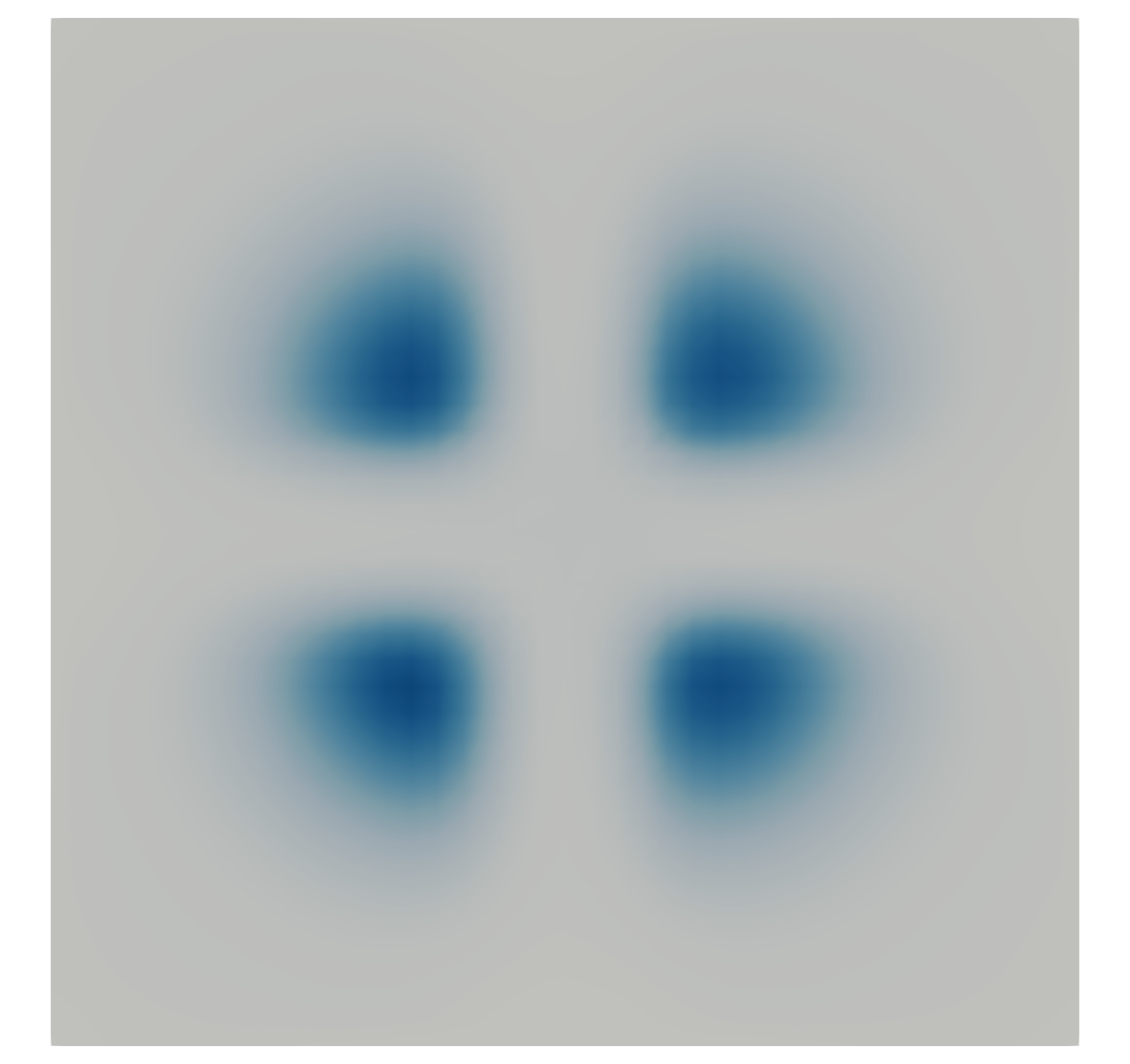}
\caption{$t=0.4$}
\end{subfigure}\\

\begin{subfigure}{0.32\textwidth}
\includegraphics[width=\textwidth]{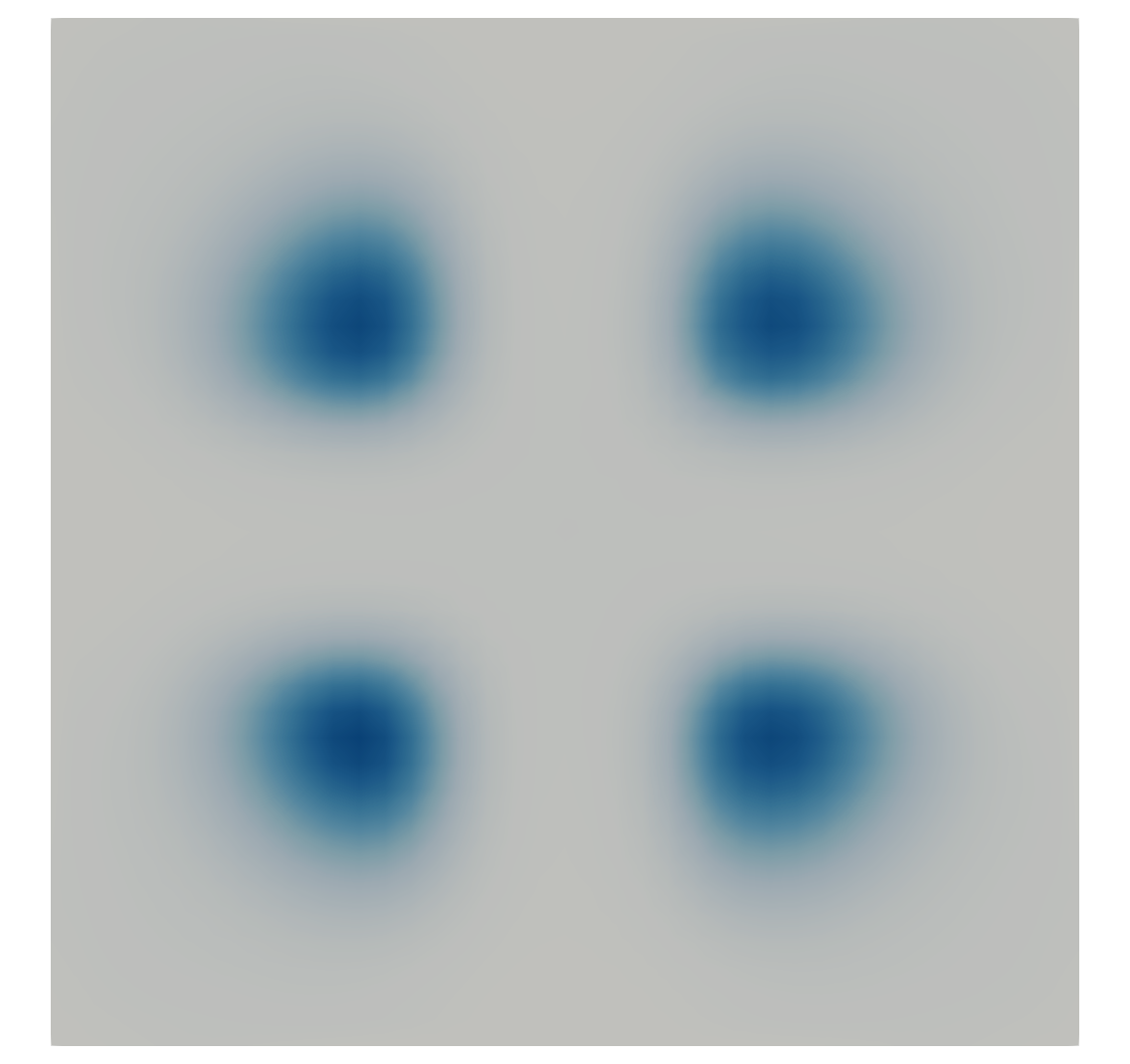}
\caption{$t=0.6$}
\end{subfigure}
\hfill
\begin{subfigure}{0.32\textwidth}
\includegraphics[width=\textwidth]{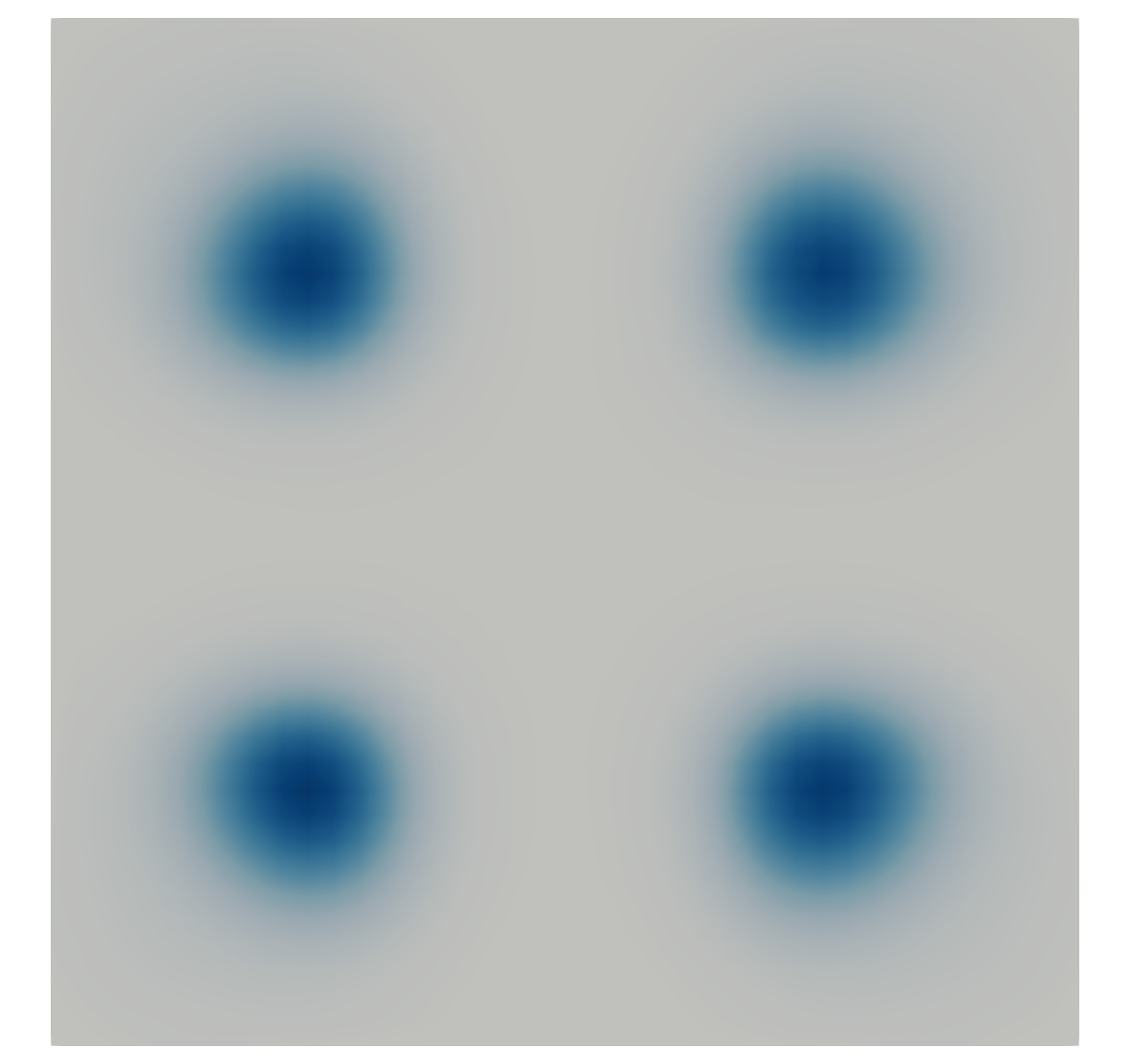}
\caption{$t=0.8$}
\end{subfigure}
\hfill
\begin{subfigure}{0.32\textwidth}
\includegraphics[width=\textwidth]{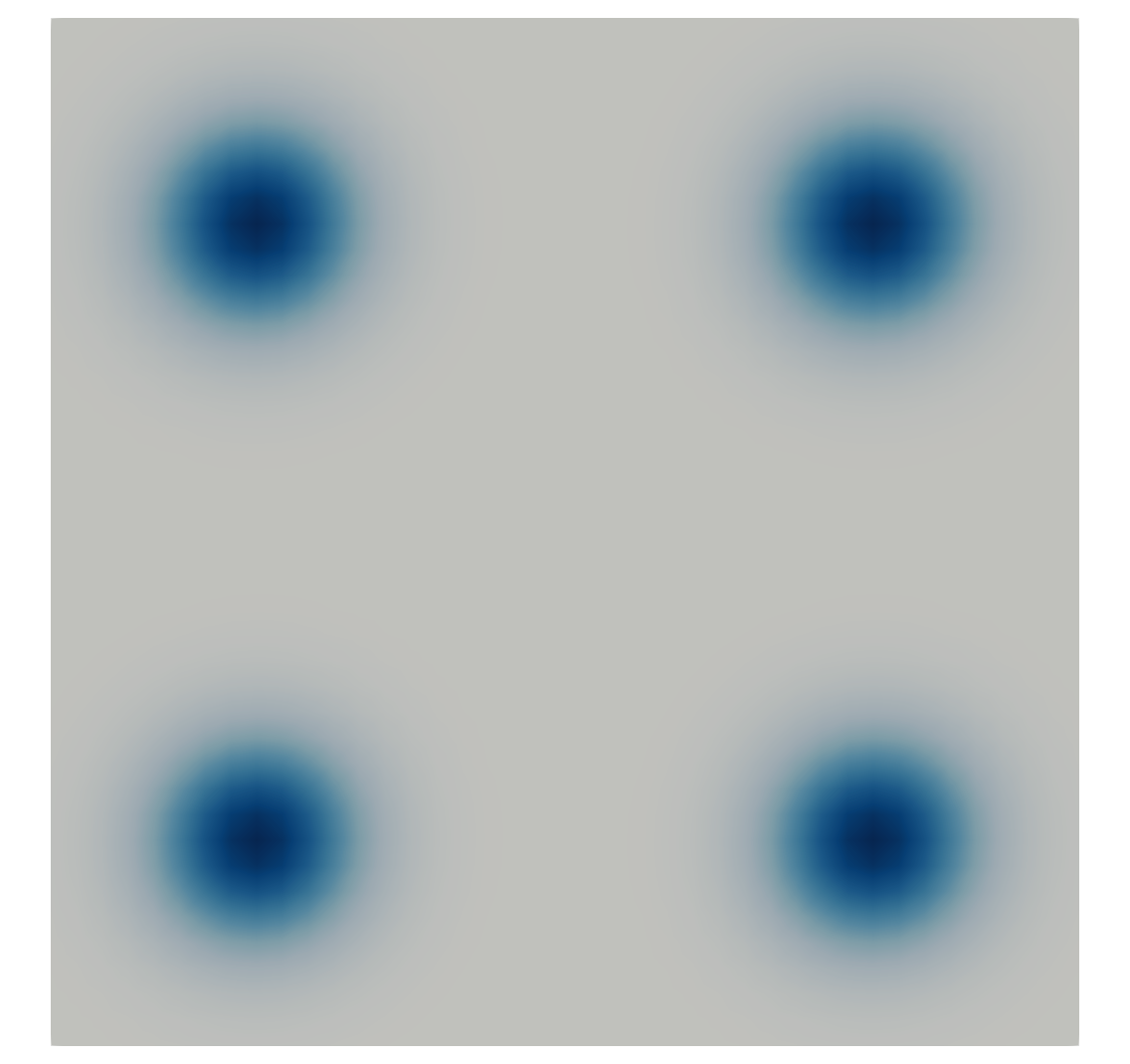}
\caption{$t=1$: $\rho_1$ distribution}
\end{subfigure}
\end{center}
\caption{Optimal transport between two distributions}
\label{optimal-transport-results}
\end{figure}

Numerical results for $\rho(\boldsymbol{x},t)$ at different time slices $t$ are represented in Figure \ref{optimal-transport-results} for $\rho_0$ being a Gaussian distribution of standard deviation 0.2 and $\rho_1$ being four identical Gaussian distributions of standard deviation 0.1 located on four opposite points of $[0;1]^2$. It can be observed how the optimal transport splits the initial distribution $\rho_0$ into four parts driving towards $\rho_1$.

\section{Conclusions and perspectives}
With the Python package \texttt{fenics\_optim} based on the FEniCS project, we propose a way to easily formulate convex variational problems arising in many applications of applied mathematics, image processing or mechanics. Convex optimization problems are formulated to fit into the conic programming framework in order to use efficient interior-point solvers especially tailored for such classes of problem. In the current form of the project, we use Mosek as the interior-point solver but other solvers could well be interfaced with the obtained discrete problems. The key point for fitting into the conic programming framework relies on a conic reformulation of convex functions. We have shown that many elementary convex functions such as $L^p$ norms arising in applications can be indeed reformulated in such a way. In the gallery of examples we tackled, we showed that various problems can be formulated with the \texttt{fenics\_optim} library in a very condensed manner. Besides, despite the fact that interior-point solvers are not necessarily the method of choice for all the considered examples, in particular for image processing applications, they are still very efficient and robust. They are therefore a good choice for a general-purpose solver for the present package. Finally, the versatility of FEniCS in terms of discretization solutions allowed to formulate very easily different discretization strategies, in particular including DG finite-elements or $H(\div)$-conforming elements which naturally arise in dual problems.\\

Obviously, there exist many aspects for improving the scope of the library or its efficiency. In particular, the current implementation of \texttt{fenics\_optim} does not allow for problems involving SDP matrix variables (Semi-Definite Programming problems), there are important applications in 3D computational mechanics which would benefit from such a feature. Besides, since version 9 of Mosek, power and exponential cones \cite{chares2009cones} are also available, which would broaden even more the class of conic representable functions including, for instance, $L^p$-norms with $p\notin\{1,2,\infty\}$,  exponential, entropy functions, etc\footnote{see \url{https://docs.mosek.com/modeling-cookbook/index.html}}. Including such a feature would therefore be a huge added value for many applications.\\
As regards computational efficiency, we mentioned that interior-point solvers, although being efficient and robust, have important memory requirements for large-scale problems since they rely on solving Newton-like systems using direct solvers. Image processing applications usually rely on proximal algorithms for solving the corresponding optimization problems, it would therefore be interesting to implement such algorithms in the package. Finally, there are some internal limitations due to the current status of the FEniCS library (e.g. Lagrange multipliers cannot be defined on one sub-part, the boundary for instance, of the domain) that could be improved. Fortunately, the FEniCS project is currently experiencing a major redevelopment to bring new functionality and improve efficiency\footnote{\url{https://github.com/FEniCS/dolfinx} and \url{https://github.com/FEniCS/ffcx/}}. We will therefore aim at taking advantage of these new developments in the later versions of the \texttt{fenics\_optim} package.

\begin{acks}
The author would like to thank Gabriel Peyré for his useful advices when preparing the manuscript and F. Bleyer for providing some of the illustrative examples input material.
\end{acks}

\appendix
\section{General principles of \texttt{fenics\_optim} internal structure}\label{appendix:structure}

\subsection{Block-structure of the problem}
The formulation of an optimization problem using \texttt{fenics\_optim} relies on a block-structure definition of variables and constraints. Let us consider, for instance formulation \eqref{obstacle-discr-block}. The conic reformulation leads to the introduction of variables $\widehat{\mathbf{y}}$ in addition to $\mathbf{u}$. Problem \eqref{obstacle-discr-block} therefore contains a block-structure of $p=2$ variables $(\mathbf{x}^1,\mathbf{x}^2)=(\mathbf{u},\widehat{\mathbf{y}})\in\RR^N\times\RR^{4M}$. The internal machinery of \texttt{fenics\_optim} works by adding sequentially new optimization variables, possibly associated with bound or conic constraints, and new linear equality or inequality constraints. Pseudo-code for defining such a block-wise structure would look like:
\begin{pythoncode}
# Problem initialization
prob = MosekProblem("My problem")

# Adding a first block variable x1
x1 = prob.add_var(V1, lx=lx1, ux=ux1, cone=K1)
# Adding a first linear constraint
prob.add_ineq_constraint(W1, A=[a1], bl=bl1, bu=bu1)
\end{pythoncode}

At this stage, the \texttt{prob} instance represents the following problem:
\begin{equation}
\begin{array}{rl}
\displaystyle{\min_{\mathbf{x}^1\in V_1}} & 0 \\
\text{s.t.} & \mathbf{l}_x^1 \leq\mathbf{x}^1 \leq \mathbf{u}_x^1 \\
& \mathbf{b}_l^1 \leq \mathbf{A}^1\mathbf{x}^1 \leq \mathbf{b}_u^1
\end{array}
\end{equation}
where $V_1$ would be $\RR^{d_1}$ in a purely discrete setting but will, in fact, be the variable \texttt{FunctionSpace} in the FEniCS FE-discretization setting. Bounds like $\mathbf{l}_x^1,\mathbf{u}_x^1, \mathbf{b}_l^1, \mathbf{b}_u^1$ are $\pm \infty$ by default (\texttt{None} in Python) and can be ignored in such case. $\Kk^1$ is a \texttt{Cone} object describing the type of cone to which the variable belongs (again \texttt{None} by default if there is no conic constraint). Finally, the linear constraint matrix $\mathbf{A}^1$ is represented by a bilinear form $a^1(y^1,x^1)$ on $W_1\times V_1$ where $y^1$ is the constraint Lagrange multiplier and $W_1$ its corresponding \texttt{FunctionSpace}. The bilinear form $a^1$ is then assembled by FEniCS to produce the discrete matrix $\mathbf{A}^1$ stored in sparse format. \\

Adding a second variable is then similar, except that constraints must now include the block-structure of both variables such as:
\begin{pythoncode}
# Problem initialization
prob = MosekProblem("My problem")

# Adding a first block variable x1
x1 = prob.add_var(V1, lx=lx1, ux=ux1, cone=K1)
# Adding a first linear constraint
prob.add_ineq_constraint(W1, A=[a1], bl=bl1, bu=bu1)

# Adding a second block variable x2
x2 = prob.add_var(V2, lx=lx2, ux=ux2, cone=K2)
# Adding a second linear constraint
prob.add_ineq_constraint(W2, A=[a21, a22], bl=bl2, bu=bu2)
\end{pythoncode}
where we now have two bilinear forms $a^{21}(y^2,x^1)$ on $W_2\times V_1$ and $a^{22}(y^2,x^2)$ on $W_2\times V_2$ leading to:
  \begin{equation}
\begin{array}{rl}
\displaystyle{\min_{(\mathbf{x}^1,\mathbf{x}^2) \in V_1\times V_2}} & 0 \\
\text{s.t.} & \mathbf{l}_x^1 \leq\mathbf{x}^1 \leq \mathbf{u}_x^1 \\
& \mathbf{b}_l^1 \leq \mathbf{A}^1\mathbf{x}^1 \leq \mathbf{b}_u^1 \\
 & \mathbf{l}_x^2 \leq\mathbf{x}^2 \leq \mathbf{u}_x^2 \\
& \mathbf{b}_l^2 \leq \mathbf{A}^{21}\mathbf{x}^1+\mathbf{A}^{22}\mathbf{x}^2 \leq \mathbf{b}_u^2
\end{array}
\end{equation}

Finally, a linear objective term $(\mathbf{c}^1)\T\mathbf{x}^1+(\mathbf{c}^2)\T\mathbf{x}^2$ can be added as 
\begin{pythoncode}
prob.add_obj_fun([c1, c2])
\end{pythoncode}

The final block-structure for a problem with $p$ blocks will therefore look like:
  \begin{equation}
\begin{array}{rl}
\displaystyle{\min_{\mathbf{x}=(\mathbf{x}^1,\ldots,\mathbf{x}^p) \in V_1\times\ldots\times V_p}} & (\mathbf{c}^1,\ldots,\mathbf{c}^p)\T(\mathbf{x}^1,\ldots,\mathbf{x}^p) \\
\text{s.t.} & \mathbf{l}_x \leq\mathbf{x} \leq \mathbf{u}_x \\
& \mathbf{b}_l \leq 
\begin{bmatrix} \mathbf{A}^{11} & 0 & \ldots & 0 \\
\mathbf{A}^{21} & \mathbf{A}^{22} & \ldots & 0\\
\vdots & \vdots & \ddots & \vdots \\
\mathbf{A}^{p1} & \mathbf{A}^{p2} & \dots & \mathbf{A}^{pp}
\end{bmatrix} \mathbf{x} \leq \mathbf{b}_u
\end{array}
\end{equation}
Note that when defining sequentially the block-constraints until variable $\mathbf{x}^i$, the blocks $\mathbf{A}^{ij}$ with $j>i$ are automatically zero since variables $\mathbf{x}^j$ have not been defined yet. This is by no means a restriction since one could perfectly define the constraint matrices once all variables have been defined. This lower-triangular structure allows however for an easier definition of the constraints in many cases. Note also that empty blocks can also be written with \texttt{0} or \texttt{None} in Python. Such symbols must be explicitly used for all $\mathbf{A}^{ij}$ with $j\leq i$.  

\subsection{Explicit construction of problem \eqref{obstacle-discr-block}}\label{appendix:obstacle}
Going back to problem \eqref{obstacle-discr-block}, one could do first:
\begin{pythoncode}
# Problem initialization
prob = MosekProblem("Obstacle problem")

# Adding a first block variable u
u = prob.add_var(Vu, lx=g)
\end{pythoncode}
creating only variable $\mathbf{u}$ and its lower bound constraint $\mathbf{u}\geq\mathbf{g}$.

Auxiliary variable $\widehat{\mathbf{y}}$ corresponds to a 4-dimensional vectorial field with degrees of freedom located at quadrature points. FEniCS provides such a functional space through the concept of \texttt{Quadrature} elements. We will use one, noted \texttt{V2}, of dimension 4 for $\widehat{\mathbf{y}}$ and one, noted \texttt{W}, of dimension 3 for the Lagrange multipliers corresponding to constraints:
\begin{align*}
y_{g,1} &= 1 \\
\mathbf{B}_g\mathbf{u} - \begin{bmatrix} y_{g,2} \\ y_{g,3} \end{bmatrix} &= 0
\end{align*}
Indeed, satisfying the above constraints for all Gauss points $g$ is equivalent to writing:
\begin{equation}
a^{12}(z,u)+a^{22}(z,y) = b(z) \quad \forall z\in W
\end{equation}
where 
\begin{align*}
a^{12}(z,u) &= \int_{\Omega} \begin{pmatrix} z_2 \\ z_3 \end{pmatrix}\cdot\nabla u \dx \\
a^{22}(z,y) &= -\int_{\Omega} z\cdot \begin{pmatrix} y_1 \\ y_2 \\ y_3 \end{pmatrix} \dx \\
b(z) &= \int_{\Omega} z \cdot \begin{pmatrix} -1 \\ 0 \\ 0 \end{pmatrix}\dx
\end{align*}
in which the integrals are computed using the same quadrature used for defining $y\in V_2$ and $z \in W$. This results in the following code:
\begin{pythoncode}
def quad_element(degree=0, dim=1):
    return VectorElement("Quadrature", mesh.ufl_cell(), 
    		degree=degree, dim=dim, quad_scheme="default")
V2 = FunctionSpace(mesh, quad_element(degree=0, dim=4))
W = FunctionSpace(mesh, quad_element(degree=0, dim=3))
y = prob.add_var(V2, cone = RQuad(4))

dxq = Measure("dx", metadata={"quadrature_scheme":"default",
							  "quadrature_degree":0})
def constraint(z):
    g = grad(u)
    a21 = dot(z, as_vector([0, g[0], g[1]]))*dxq
    a22 = -dot(z, as_vector([y[1], y[2], y[3]]))*dxq
    return [a21, a22]
def rhs(z):
    return -z[0]*dxq
prob.add_eq_constraint(W, A=constraint, b=rhs)
\end{pythoncode}
where $y\in V_2$ is created by specifying that it also belongs to a rotated quadratic cone $\Qq_r^4$. This statement is understood point-wise, meaning that at each degree of freedom (Gauss point) location $x_g$, the local 4-d vector $y(x_g)$ belongs to $\Qq_r^4$.
The \texttt{dxq} measure is used to enforce a one-point quadrature on each cell, the same used for the definition of \texttt{V2} and \texttt{W}. Finally, the constraint matrix is passed as a function (\texttt{constraint}) of the Lagrange multiplier $z\in W$ and returns a list of 2 bilinear forms corresponding to both blocks in $u$ and $y$, while the constraint right-hand side is also passed as a function of $z$ (\texttt{rhs}) and returns a single linear form in $z$. A similar syntax would be used for inequality constraints.

Finally, the objective term is set as a list of two linear forms in $u$ and $y$ respectively:
\begin{pythoncode}
prob.add_obj_func([-dot(load, u)*dx, y[0]*dxq])
\end{pythoncode}
Note again the use of the one-point quadrature measure for the second term. 

One role of \texttt{ConvexFunction} classes described in \ref{sec:fenics-obstacle} is to avoid for the user to explicitly define
function spaces for the additional variables and Lagrange multipliers. 
The complete script of this implementation can be found in \texttt{demos/obstacle/obstacle\_problem\_explicit\_construction}.

\section{Conic reformulation of problem \eqref{BENDING-PLATE}}\label{appendix:plate}
We consider function $\pi: M\in\mathbb{S}^2 \longleftrightarrow \frac{2m}{\sqrt{3}}\sqrt{M_{11}^2+M_{22}^2+M_{12}^2+M_{11}M_{22}}$. Expressing $M\in \mathbb{S}^2$ as $\mathbf{X} = (M_{11}, M_{22}, 2M_{12})$, we have that:
\begin{equation}
\pi(M) = \widehat{\pi}(\mathbf{X})=\frac{m}{\sqrt{3}}\sqrt{\mathbf{X}\T\mathbf{C}\mathbf{X}} \text{ with } \mathbf{C} = \begin{bmatrix}
4 & 2 & 0 \\ 2 & 4 & 0 \\ 0 & 0 & 1 
\end{bmatrix}
\end{equation}
Computing the Cholesky factor $\mathbf{J}= \begin{bmatrix} 2 & 1 & 0 \\ 0 & \sqrt{3} & 0 \\ 0 & 0 & 1 \end{bmatrix}$  of matrix $\mathbf{C}$, we have that $\widehat{\pi}(\mathbf{X}) = \dfrac{m}{\sqrt{3}}\sqrt{\mathbf{X}\T\mathbf{J}\T\mathbf{J}\mathbf{X}} = \dfrac{m}{\sqrt{3}}\|\mathbf{JX}\|_2$.

\bibliographystyle{apalike}
\bibliography{fenics_optim}

\end{document}